\documentclass[11pt,reqno]{amsart}
\usepackage{epsfig}
\usepackage{amsthm}
\usepackage{amsmath}
\usepackage{amssymb}
\usepackage{amsbsy}
\usepackage{amsfonts}
\usepackage{latexsym}
\usepackage{indentfirst}
\usepackage{textcomp}
\newtheorem{thm}{Theorem}[section]
\newtheorem{cor}[thm]{Corollary}
\newtheorem{lem}[thm]{Lemma}

\theoremstyle{definition}
\newtheorem{defn}[thm]{Definition}
\theoremstyle{remark}
\newtheorem{rem}[thm]{Remark}
\numberwithin{equation}{section}

\makeatletter
\renewcommand{\@evenhead}{\thepage \hfil V.V.BORZOV \hfil}
\renewcommand{\@oddhead}{\hfil ORTHOGONAL POLYNOMIALS AND  OSCILLATOR  ALGEBRAS \hfil \thepage}
\begin{document}
\begin{center}

{\Large \bf ORTHOGONAL POLYNOMIALS AND  GENERALIZED \\
         OSCILLATOR  ALGEBRAS}

\vspace{5mm}

{\large \bf V.V. Borzov}

\medskip

{\em Department of Mathematics, St.Petersburg University of Telecommunications, \\
191065, Moika  61,  St.Petersburg, Russia}

\medskip

{\em (Received December 6, 1999)}

\end{center}
\bigskip
For any orthogonal polynomials system on real line we construct an appropriate oscillator algebra such that
 the polynomials make up the eigenfunctions system  of the oscillator hamiltonian. The general scheme is
divided into two types: a symmetric scheme and a non-symmetric scheme. The general approach is illustrated
by the examples of the classical orthogonal polynomials: Hermite, Jacobi and Laguerre polynomials. For these
polynomials we obtain the explicit form of the hamiltonians, the energy levels and the explicit form of the
impulse operators.

\bigskip


\medskip


\bigskip

\section{Introduction}

The connection of orthogonal polynomials with the classical groups (\cite{vil}) as well as with the quantum ones
(\cite{kor}) is well known. We discuss here the connection of orthogonal polynomials with the Heisenberg algebra
of generalized (deformed  (\cite{bie, mac, dam}) as an example) oscillator.
Recall  that (see, for example, \cite{lan}) the Hermite polynomials (after multiplication by $ \ exp(-x^2)$)
make up the eigenfunctions system of the energy operator for the quantum mechanical harmonic
oscillator. Many of the known $q$-Hermite polynomials (\cite{asw, asi, flo}) are also the eigenfunctions of the energy
operator for a deformed oscillator. It is well known that  orthogonal polynomials ,which in a sense generalize
the Hermite polynomials , appear in the analysis of  the irreducible representations of  the algebra of an appropriate
oscillator. In this paper we propose another way of looking at  the  connection of orthogonal polynomials with some
generalized oscillator algebras. Namely, given an orthogonal polynomials system, we construct an appropriate
oscillator algebra so that  the polynomials make up a eigenfunctions system  of the oscillator hamiltonian.

The aim of this paper is to present the classical orthogonal polynomials as eigenfunctions of an energy operator
 for a generalized oscillator. Let us take a brief look at the considered approach. A preassigned Hilbert space with an
 orthogonal polynomials systems  (for instance, one of the above-mentioned classical polynomials systems)
as a basis is considered as a Fock space.  As it usually is, we define the ladder operators (annihilation) $a^-$
and (creation) $a^+$ as well as the number operator $N $ in this space. By a standard manner we use these
operators to build up the following selfadjoint operators: the position operator   $X$, the momentum operator
$P$ as well as the energy operator (hamiltonian) $H={X^2}+{P^2}$. By analogy with the usual Heisenberg algebra
these operators generate an algebra, which naturally is called a generalized oscillator algebra.
It turns out that  the operator $H$ has a simple discrete spectrum. The initial orthogonal polynomials set
is an eigenfunctions system of the energy operator $H$. Via the Poisson kernel of this system is determined
a generalized Fourier transform, which establishes the usual link between the operators $X$ and $P$.
The energy operator $H$ is invariable under the action of this transform.
The explicit form of the Poisson kernels for  the classical orthogonal polynomials (the analog of the Mehler
 formula \cite{meh}) see in \cite{wat1}, \cite{wat2}, \cite{wat3}.
The orthogonal polynomials systems (OPS) can be further divided into two types: symmetric systems and
non-symmetric systems. OPS is called a symmetric system if the orthogonality measure for these polynomials
is symmetric about the origin; otherwise it is called a non-symmetric system. In the former case the Jacobi
matrix of the operator   $X$ (in the Fock representation) has the trivial diagonal. Note that the above-mentioned
oscillator algebra arise only in the first case. In the latter case one can also construct a generalized oscillator
algebra. However the oscillator hamiltonian takes the standard form only in new "coordinate-impulse" operators
, which can result from the previous operators $X$ and $P$ by a rotation.

\section{Symmetric scheme}

\subsection{}  Let $\mu$ be a positive Borel measure on the real line $R^1$ such that
\begin{equation}
\int_{-\infty}^{\infty}{\mu(dx)}=1,\qquad
\mu_{2k+1}=\int_{-\infty}^{\infty}{x^{2k+1}\mu(dx)}=0,\qquad
k= 0,1,\dots  .
\label{triv1}
\end{equation}
The measure $\mu$ is called a symmetric probability measure. By $\tt H$ we denote the Hilbert space $L^2(R^1;{\mu)}.$
Let $\left\{ {b_{n}}\right\}_{n=0}^{\infty }$,\quad ${b_{n}>0}$,\quad ${n=0,1,\dots }$ be a positive sequence defined by
the algebraic equations system
\begin{equation}
\mu_{2k}= {{b _{0}^2}\cdot( {b _{0}^2}+{b _{1}^2}) \cdots \sum_{j=0}^{k-1}{{b _{j}^2}} },\qquad
\label{triv3}
k= 0,1,\dots ,
\end{equation}
where
\begin{equation}
\mu_{0}=1,\qquad
\mu_{2k}=\int_{-\infty}^{\infty}{x^{2k}\mu(dx)},\qquad
k= 0,1,\dots ,
\label{triv2}
\end{equation}
Obviously, there is the unique solution to the system (\ref{triv3})
\begin{equation}
{b _{0}^2}=\mu_{2},\qquad
{b _{1}^2}=\frac{\mu_{4}}{\mu_{2}}-\mu_{2},\dots .
\label{triv4}
\end{equation}

\begin{defn}
A polynomial set $\left\{ {\psi_{n}(x)}\right\}_{n=0}^{\infty }$
is called a canonical  polynomial system  if it is defined by the following recurrence relations:

\begin{equation}
{x {\psi_{n}(x)}}= {b_{n}{\psi_{n+1}(x)}}+{b_{n-1}{\psi_{n-1}(x)}},\qquad  n\geq 0,\qquad{ b_{-1}}=0,
\label{triv5}
\end{equation}

\begin{equation}
{\psi_{0}(x)}=1,
\label{triv6}
\end{equation}
where the positive sequence $\left\{ {b_{n}}\right\}_{n=0}^{\infty }$ is given.

\end{defn}

\begin{rem}
\label{j}
1. The canonical  polynomial system $\left\{ {\psi_{n}(x)}\right\}_{n=0}^{\infty }$
is uniquely determined by the symmetric probability measure  $\mu$.

2.The recurrence relations  (\ref{triv5}) give us the symmetric Jacobi matrix
$$
J=\left\{
b_{ij}
\right\}
_{i,j=0}^{\infty }
$$
which has the positive elements $b_{i,i+1}=b_{i+1,i},\quad i=0,1,$\dots only distinct from zero.
If the moment problem (\cite{ach}) for the matrix $J$ is a determined one, then the canonical  polynomial system
$\left\{ {\psi_{n}(x)}\right\}_{n=0}^{\infty }$ is completed in the space $\tt H$. Otherwise (when the moment problem
is a undetermined one) the canonical  polynomial system $\left\{ {\psi_{n}(x)}\right\}_{n=0}^{\infty }$
is completed in the space $\tt H$ if and only if the measure  $\mu$ is a $N$- extremal solution (\cite{ach})
of the moment problem for the matrix $J$.
\end{rem}

The following theorem is true.
\begin{thm}\label{a}
Let
$\left\{
{\psi_{n}(x)}
\right\}_{n=0}^{\infty }$
be a set of real polynomials
satisfying   recurrence relations  (\ref{triv5}) and a initial condition
(\ref{triv6}); let $\mu$ be a symmetric probability measure on the real
line $R^1$, that is the conditions (\ref{triv1})
for $\mu$ are valid. The set
$\left\{
{\psi_{n}(x)}
\right\}_{n=0}^{\infty }$ is a system of polynomials orthonormal
with respect to the measure
$\mu$ if and only if a positive sequence
$\left\{
{b_{n}}\right\}_{n=0}^{\infty }$ involved
in the recurrence relations  (\ref{triv5}) is a solution of the system
(\ref{triv3}), where $\mu_{2k}$ are defined by
(\ref{triv2}).
\end{thm}

\subsection{}

Let $\psi(x)$ be a real-valued function such that $\frac{1}{\psi(x)}$ is  measurable with respect to the
above-mentioned measure  $\mu$. Let us introduce  new measure $\nu$ realised by
\begin{equation}
{\nu(dx)}={\left|{\psi(x)}\right|^{-2}}{\mu(dx)}.
\label{triv7}
\end{equation}

\begin{rem}
1. The function
${\left| \psi(x) \right|}^{-2}$ is locally integrable but it is not necessarily that
${\left| \psi(x) \right|}^{-2}\in{L^1(R^1;{\mu(dx)}}$.

2.In general, the conditions (\ref{triv1}) break down for the measure  $\nu$.
\end{rem}
Now we consider  another Hilbert space \begin{math}{\tt G}=L^2(R^1;{\nu(dx)})\end{math} with the measure  $\nu$
defined by (\ref{triv7}). We define the functions system
$\left\{ {\phi_{n}(x)}\right\}_{n=0}^{\infty }$, ${\phi_{n}(x)}\in{\tt G},n= 0,1,\dots$ by
\begin{equation}
{\phi_{n}(x)}={\psi(x)}{\psi_{n}(x)},\qquad
\label{triv8}
n= 0,1,\dots ,
\end{equation}
where the set $\left\{ {\psi_{n}(x)}\right\}_{n=0}^{\infty }$ is a canonical  polynomial system in above space $\tt H$.

The following statement is a simple consequence of the theorem {\ref{a}}.
\begin{cor}
\label{ded}
If the system $\left\{ {\psi_{n}(x)}\right\}_{n=0}^{\infty }$ is a canonical system of polynomials orthonormal
with respect to the measure  $\mu$ in the space $\tt H$, then the set
$\left\{ {\phi_{n}(x)}\right\}_{n=0}^{\infty }$ ,${\phi_{n}(x)}\in{\tt G},n= 0,1,\dots$ defined by (\ref{triv8}) is a
orthonormal system in the space \begin{math}{\tt G}=L^2(R^1;{\nu(dx))}\end{math}. Besides, this system satisfies
the same recurrence relations  (\ref{triv5}) and the initial condition
\begin{equation}
{\phi_{0}(x)}={\psi(x)}.
\label{triv9}
\end{equation}
\end{cor}

\begin{rem}
It is evident, that a completeness of system $\left\{ {\phi_{n}(x)}\right\}_{n=0}^{\infty }$ in the space $\tt G$ is
equivalent to the one of $\left\{ {\psi_{n}(x)}\right\}_{n=0}^{\infty }$ in the space $\tt H$.
\end{rem}

\subsection{The Poisson kernel}

For the reader's convenience we remind the definition of the Poisson kernel in the Hilbert space
\begin{math}{\tt F}=L^2(R^1;\rho(dx))\end{math}, where $\rho$ is a positive Borel measure on the real line $R^1$.
Let a set $\left\{ {\varphi_{n}(x)}\right\}_{n=0}^{\infty }$ be an orthonormal basis in the space ${\tt F}$. From now on
we will use the notation $\tt F$ instead of $\tt G$ or $\tt H$ if  both spaces are regarded together. Let us denote
by ${\tt F}_{1}$, ${\tt F}_{2}$ the first and second copies of the space $\tt F$ respectively:
\begin{equation}
{{\tt F}_{1}}={L^2(R^1;{\rho(dx)})},\qquad
{{\tt F}_{2}}={L^2(R^1;{\rho(dy)})}.
\label{triv10}
\end{equation}

The Poisson kernel $\mathfrak{K}_{\tt F}(x,y;t)$  on  ${\tt F}_{1}\otimes {\tt F}_{2}$ is defined by the formula
\begin{equation}
{\mathfrak{K}_{\tt F}(x,y;t)}=\sum_{n=0}^{\infty }{{t^{n}}\cdot {\varphi_{n}(x)}\cdot {\varphi_{n}(y)}}.
\label{triv11}
\end{equation}
From (\ref{triv8}) and (\ref{triv11}) it follows that:
\begin{equation}
{\mathfrak{K}_{\tt G}(x,y;t)}={\psi(x)}{\psi(y)}{\mathfrak{K}_{\tt H}(x,y;t)}.
\label{triv12}
\end{equation}
We define the integral operators $K_{\tt F}$: ${\tt F}_{1} \longmapsto{\tt F}_{2}$ and
${K^{\prime}} _{\tt F}$ : ${\tt F}_{2} \longmapsto{\tt F}_{1}$ by the following formulas
\begin{equation}
(K_{\tt F}{f})(y)=\int_{-\infty}^{\infty}{f(x)\mathfrak{K}_{\tt F}(x,y;t)\rho(dx)},
\label{triv13}
\end{equation}
\begin{equation}
({K^{\prime}}_{\tt F}{g})(x)=\int_{-\infty}^{\infty}{g(y)\overline{\mathfrak{K}_{\tt F}(x,y;t)}\rho(dy)}.
\label{triv14}
\end{equation}
It is easy to prove the following lemmas.
\begin{lem}
Let the set $\left\{ {\varphi_{n}(x)}\right\}_{n=0}^{\infty }$ be an orthonormal system in the space ${\tt F}_{1}$.
If this system is completed in ${\tt F}_{1}$ and $|t|=1$, then the integral operators (\ref{triv13}),  (\ref{triv14})
are unitary ones:
\begin{equation}
{K^{\prime}}_{\tt F}={K^{*}}_{\tt F}={K^{-1}}_{\tt F}.
\label{triv15}
\end{equation}
\end{lem}

\defn
The operator $U_x$: ${\tt H}_{1} \longmapsto{\tt G}_{1}$ is defined by
\begin{equation}
f(x)={U_x}e(x)\Leftrightarrow
f(x)={\psi(x)}e(x),\qquad
e(x) \in {\tt H}_{1},
f(x) \in {\tt G}_{1}.
\label{triv16}
\end{equation}
Likewise, the operator $U_y$: ${\tt H}_{2} \longmapsto{\tt G}_{2}$ is determined.

\begin{lem}
If the set $\left\{ {\psi_{n}(x)}\right\}_{n=0}^{\infty }$ is an orthonormal basis in the space ${\tt H}_{1}$,
then the set
$\left\{ {\phi_{n}(x)}\right\}_{n=0}^{\infty }$ , defined by (\ref{triv8}) is an
orthonormal basis in the space ${\tt G}_{1}$. Besides, the operator $U_x$ determined by (\ref{triv16}) is an unitary one
\begin{equation}
{U_x}^{*}={U_x}^{-1}.
\label{triv17}
\end{equation}
The same affirmation is true for the operator $U_y$.
\label{c}
\end{lem}

\begin{lem}
The operators  (\ref{triv13}) and (\ref{triv16}) satisfy the following relations:
\begin{equation}
K_{\tt G}={U_y}{K_{\tt H}}{U_x}^{-1},\qquad
K_{\tt H}={U_y}^{-1}{K_{\tt G}}{U_x}.
\label{triv18}
\end{equation}
\label{u}
\end{lem}
\begin{proof}
The proof is trivial.
\end{proof}

\subsection{The Hamiltonian formulation}

From  now on we assume that the orthonormal system $\left\{ {\varphi_{n}(x)}\right\}_{n=0}^{\infty }$ is completed in
${\tt F}_{1}=L^{2}(R^1;{\rho (dx)})$. The relations (\ref{triv5}) indicate a manner by which the position operator
$X_{{\tt F}_{1}}$ acts on the elements of this basis in the Fock space ${\tt F}_{1}$. Let us remember (\cite{bir})
that the domain $D(X_{{\tt F}_{1}})$  of operator $X_{{\tt F}_{1}}$ is defined by
\begin{equation}
D(X_{{\tt F}_{1}})=\left\{ f(x)\in{\tt F}_{1}|{\int_{-\infty}^{\infty}{{|f(x)|^{2}}({1+x^{2}}){\rho(dx)}}<{ \infty}}\right\}
\label{triv19}
\end{equation}
Using  (\ref{triv13}), (\ref{triv14}), we define now a momentum operator $P_{{\tt F}_{1}}$,
which is conjugate to the position operator $X_{{\tt F}_{1}}$ with respect to the basis
$\left\{ {\varphi_{n}(x)}\right\}_{n=0}^{\infty }$ of  ${\tt F}_{1}$ in the following way:
\begin{equation}
P_{{\tt F}_{1}}={K_{\tt F}}^{*}{Y_{{\tt F}_{2}}}{K_{\tt F}}.
\label{triv20}
\end{equation}
Note that a operator $Y_{{\tt F}_{2}}$ in  (\ref{triv20}) is a position operator in the space ${\tt F}_{2}$
 defined by analogy with the formulas (\ref{triv5}). In general , we have ($|t|=1$)
\begin{equation}
D(P_{{\tt F}_{1}})={K_{\tt F}}^{*}D(Y_{{\tt F}_{2}}).
\label{triv21}
\end{equation}
Finally, we define the operator
\begin{equation}
{H_{{\tt F}_{1}}}(t)=(X_{{\tt F}_{1}})^{2}+({P_{{\tt F}_{1}}}(t))^{2}.
\label{triv22}
\end{equation}
The following theorem is our main result of  the present section. The proof  is very simple and it is omitted.
\begin{thm}
Let a canonical  (polynomial) system $\left\{ {\varphi_{n}(x)}\right\}_{n=0}^{\infty }$ be completed in
the space ${\tt F}_{1}$. This system is a set of  eigenfunctions of the selfadjoint operator
${H_{{\tt F}_{1}}}(t)$ in ${\tt F}_{1}$ defined by (\ref{triv22}) in ${\tt F}_{1}$   if and only if  $t=\pm\imath$.
Moreover, the eigenvalues of the operators ${H_{{\tt F}_{1}}}(\pm\imath)$ are equal to
\begin{equation}
\lambda_{0}=2{b_{0}^{2}},\qquad
\lambda_{n}=2({b_{n-1}^{2}}+{b_{n}^{2}}),\qquad
n\geq 1.
\label{triv23}
\end{equation}
\label{b}
\end{thm}

\begin{rem}
The operator ${H_{{\tt F}_{1}}}={H_{{\tt F}_{1}}}(-\imath)$ is said to be a hamiltonian of the orthonormal system
$\left\{ {\varphi_{n}(x)}\right\}_{n=0}^{\infty }$. The domain of the operator $H_{{\tt F}_{1}}$ is obtained from
(\ref{triv19}) and (\ref{triv21}) by the following formulas
$$
D(H_{{\tt F}_{1}})=
\overline{D(X_{{\tt F}_{1}}^{2})
\cap
D(P_{{\tt F}_{1}}^{2})}.
$$
Also, we denote by
\begin{equation}
P_{{\tt F}_{1}}={P_{{\tt F}_{1}}}(-\imath),\quad
P_{{\tt F}_{2}}={P_{{\tt F}_{2}}}(\imath),\quad
H_{{\tt F}_{2}}={H_{{\tt F}_{2}}}(\imath).
\label{triv24}
\end{equation}
\end{rem}
The proof of the following lemmas is left to the reader.
\begin{lem}
\label{du}
The operators $X_{{\tt F}_{1}}$, $P_{{\tt F}_{1}}$,$H_{{\tt F}_{1}}$
act on the basis vectors $\left\{ {\varphi_{n}(x)}\right\}_{n=0}^{\infty }$ of the space ${\tt F}_{1}$ by
\begin{align}
\label{triv25}
X_{{\tt F}_{1}}{\varphi_{0}(x)}&={b_{0}}{\varphi_{1}(x)},\\
\label{triv26}
P_{{\tt F}_{1}}{\varphi_{0}(x)}&=-{\imath}{b_{0}}{\varphi_{1}(x)},\\
\label{triv27}
H_{{\tt F}_{1}}{\varphi_{0}(x)}&={\lambda_{0}}{\varphi_{0}(x)},\\
\label{triv28}
X_{{\tt F}_{1}}{\varphi_{n}(x)}&={b_{n-1}}{\varphi_{n-1}(x)}+{b_{n}}{\varphi_{n+1}(x)},\quad
n\geq 1,\\
\label{triv29}
P_{{\tt F}_{1}}{\varphi_{n}(x)}&={\imath}({b_{n-1}\varphi_{n-1}(x)}-{{b_{n}\varphi_{n+1}(x))}},\quad
n\geq 1,\\
\label{triv30}
H_{{\tt F}_{1}}{\varphi_{n}(x)}&={\lambda_{n}}{\varphi_{n}(x)},\quad
\end{align}
where the eigenvalues $\lambda_{n}, n\geq  0,$ are defined by (\ref{triv23}).
\end{lem}
\begin{lem}
Under the assumptions of the lemma {\ref{c}} the operators (\ref{triv22}),(\ref{triv20}) comply with
the following relations:
\begin{equation}
P_{{\tt G}_{1}}={U_x}{P_{{\tt H}_{1}}}{U_x}^{-1},\qquad
H_{{\tt G}_{1}}={U_x}{H_{{\tt H}_{1}}}{U_x}^{-1}.
\label{triv31}
\end{equation}
\end{lem}
\begin{rem}
The previous statement  still stands for the operators $(P_{{\tt F}_{1}})(t), (H_{{\tt F}_{1}}))(t)$  at
any  $t$ ($|t|=1$).
\end{rem}

\subsection{The generalised Fourier transform}

In this subsection  we define the Fourier transform conforming to an orthonormal system
$\left\{ {\varphi_{n}(x)}\right\}_{n=0}^{\infty }$ in the space ${\tt F}_{1}$(see{\cite{ars1}}).
\defn  Let    $\left\{ {\varphi_{n}(x)}\right\}_{n=0}^{\infty }$ be an orthonormal basis in the space ${\tt F}_{1}$.
The unitary operators ${K_{\tt F}}(\pm\imath)$ are called the generalized (direct and inverse) Fourier transforms.
We denote by
\begin{equation}
F_{\varphi}={K_{\tt F}}(-\imath),\qquad
{F_{\varphi}}^{-1}={K_{\tt F}}(\imath).
\label{triv32}
\end{equation}
The following theorem can be proved by direct calculations.
\begin{thm}
We have in the Hilbert space   ${\tt F}_{2}$
for the operators (\ref{triv20}):
\begin{equation}
P_{{\tt F}_{2}}{F_{\varphi}}={F_{\varphi}}{X_{{\tt F}_{1}}},\qquad
Y_{{\tt F}_{2}}{F_{\varphi}}={F_{\varphi}}{P_{{\tt F}_{1}}},\qquad
P_{{\tt F}_{1}}{{F_{\varphi}}^{-1}}={{F_{\varphi}}^{-1}}{Y_{{\tt F}_{2}}},\qquad
X_{{\tt F}_{1}}{{F_{\varphi}}^{-1}}={{F_{\varphi}}^{-1}}{P_{{\tt F}_{2}}}.
\label{triv33}
\end{equation}
and for the operators (\ref{triv22}):
\begin{equation}
H_{{\tt F}_{2}}{F_{\varphi}}={F_{\varphi}}{H_{{\tt F}_{1}}},\qquad
H_{{\tt F}_{1}}{{F_{\varphi}}^{-1}}={{F_{\varphi}}^{-1}}{H_{{\tt F}_{2}}}.
\label{triv34}
\end{equation}
\label{d}
\end{thm}

\subsection{The generalized oscillators algebra}

Let $\left\{ {\varphi_{n}(x)}\right\}_{n=0}^{\infty }$ be an orthonormal basis in the Fock space ${\tt F}_{1}$.
We construct some (generalized)  oscillators algebra corresponding the system
$\left\{ {\varphi_{n}(x)}\right\}_{n=0}^{\infty }$. To this end we define  ladder operators
${a^{+}}_{{\tt F}_{1}}$ and ${a^{-}}_{{\tt F}_{1}}$ by the usual formulas:
\begin{equation}
a^{+}_{{\tt F}_{1}}={\frac{1}{\sqrt{2}}}
\left (
{X_{{\tt F}_{1}}}+\imath{P_{{\tt F}_{1}}}
\right )
,\qquad
a^{-}_{{\tt F}_{1}}={\frac{1}{\sqrt{2}}}
\left (
{X_{{\tt F}_{1}}}-\imath{P_{{\tt F}_{1}}}
\right ).
\label{triv35}
\end{equation}
It is readily seen that (for the classical orthogonal polynomials)
$$
{a^{-}_{{{\tt F}_{1}}}}^{*}=a^{+}_{{\tt F}_{1}},\qquad
{a^{+}_{{{\tt F}_{1}}}}^{*}={a^{-}}_{{{\tt F}_{1}}}.
$$
and
$$
D(a^{-}_{{\tt F}_{1}})=D(a^{+}_{{\tt F}_{1}})=\overline{D(X_{{\tt F}_{1}})
\cap
D(P_{{\tt F}_{1}})}.
$$
\begin{lem}
The action of operators (\ref{triv35}) on the vectors of the basis in the space ${{\tt F}_{1}}$  is given by
the standard formulas:
\begin{equation}
{a^{+}_{{\tt F}_{1}}}{\varphi_{n}(x)}={\sqrt{2}}{b_{n}}{\varphi_{n+1}(x)},\qquad
{a^{-}_{{\tt F}_{1}}}{\varphi_{n}(x)}={\sqrt{2}}{b_{n-1}}{\varphi_{n-1}(x)},\qquad
n\geq  0.
\label{triv36}
\end{equation}
\end{lem}
It is easy to prove from (\ref{triv35}), (\ref{triv29}) and (\ref{triv28}).
\begin{lem}
Under the assumptions of the lemma {\ref{c}} the operators (\ref{triv36}),(\ref{triv16}) comply with
the following relations:
\begin{equation}
{a^{\pm}_{{\tt G}_{1}}}={U_x}{a^{\pm}_{{\tt H}_{1}}}{U_x}^{-1},\qquad
[{X_{{\tt F}_{1}}},{P_{{\tt F}_{1}}}]=\imath[{a^{-}_{{\tt H}_{1}}},{a^{+}_{{\tt H}_{1}}}].
\label{triv37}
\end{equation}
\end{lem}
\defn
An operator $N_{{\tt F}_{1}}$ in the Fock space ${\tt F}_{1}$  equipped with the basis
$\left\{ {\varphi_{n}(x)}\right\}_{n=0}^{\infty }$ is called a number operator if it acts on  basis vectors by formulas:
\begin{equation}
{N_{{\tt F}_{1}}}{\varphi_{n}(x)}={n}{\varphi_{n}(x)},\qquad
n\geq  0.
\label{triv38}
\end{equation}
\begin{lem}
Under the assumptions of the lemma {\ref{c}} the operators (\ref{triv38}) satisfy the following relations:
\begin{equation}
N_{\tt{ G}_{1}}={U_x}{N_{{\tt H}_{1}}}{U_x}^{-1}.
\label{triv39}
\end{equation}
\label{f}
\end{lem}
The proof is simple.
\begin{rem}
1. We denote by $B(N)$ a function of operator $N$ in the space $\tt{ F}_{1}$
which acts on the vectors of the basis $\left\{ {\varphi_{n}(x)}\right\}_{n=0}^{\infty }$ by
\begin{equation}
B({N_{{\tt F}_{1}}}){\varphi_{n}(x)}={b_{n-1}^{2}}{\varphi_{n}(x)},\qquad
n\geq  0,\quad
b_{-1}=0.
\label{triv40}
\end{equation}
2. Let the assumptions of the lemma {\ref{c}} be held.Then  from (\ref{triv40}) and the lemma {\ref{f}} it follows that
\begin{equation}
B({N_{{\tt G}_{1}}})={U_x}B({N_{{\tt H}_{1}}}){U_x}^{-1}.
\label{triv41}
\end{equation}
\end{rem}
The following theorem is our main result of  the present subsection. The proof  is very simple and it is omitted.
\begin{thm}
Under the assumptions of the lemma {\ref{c}} the operators (\ref{triv35}), (\ref{triv40}) in the Fock space $\tt{ F}_{1}$
 satisfy the following relations:
\begin{equation}
[{{a^{-}}_{{\tt F}_{1}}},{{a^{+}}_{{\tt F}_{1}}}]=2(B(N_{{\tt F}_{1}}+I_{{\tt F}_{1}})-B(N_{{\tt F}_{1}})),\quad
[N_{{\tt F}_{1}},{{a^{\pm}}_{{\tt F}_{1}}}]=\pm{a^{\pm}_{{\tt F}_{1}}}.
\label{triv42}
\end{equation}
Let the sequence $\left\{ b_{n} \right\}_{n=0}^{\infty }$ be defined by (\ref{triv3}) in the space $\tt{ F}_{1}$
with the measure $\rho$. If there is a real number $A$ and a real function $C(n)$,
such that this sequence satisfies the following recurrence relation:
\begin{equation}
{b_{n}^{2}}-A{b_{n-1}^{2}}=C(n),\qquad
n\geq  0,\quad
b_{-1}=0,
\label{triv43}
\end{equation}
then  the operators (\ref{triv35}),(\ref{triv40}) satisfy the following conditions:
\begin{equation}
{a^{-}_{{\tt F}_{1}}}{a^{+}_{{\tt F}_{1}}}-A{a^{+}_{{\tt F}_{1}}}{a^{-}_{{\tt F}_{1}}}=2C(N_{{\tt F}_{1}}),
\label{triv44}
\end{equation}
apart from(\ref{triv42}).
Here the function $C(N)$ is defined similarly  (\ref{triv40})  with $C(n)$ instead of
${b_{n-1}}^{2}$.
\label{g}
\end{thm}
\begin{proof}
It is follows from the obvious relations:
\begin{eqnarray}
a^-_{{\tt F}_{1}}a^+_{{\tt F}_{1}}\varphi_{n}(x)=2b_n^2 \varphi_{n}(x), \nonumber\\
a^+_{{\tt F}_{1}}a^-_{{\tt F}_{1}}\varphi_{n}(x)=2b_{n-1}^2\varphi_{n}(x),\nonumber\\
n\geq  0,\quad
b_{-1}=0.
\label{triv45}
\end{eqnarray}
\end{proof}
\defn
An algebra $A_{\varphi }$  is called a generalized oscillator algebra corresponding to the orthonormal
system $\left\{ {\varphi_{n}(x)}\right\}_{n=0}^{\infty }$ if $A_{\varphi }$ is generated by  generators
$a^{\pm}_{{\tt F}_{1}}$, $N_{{\tt F}_{1}}$, which satisfy the relations of (\ref{triv45}) and the two latter ones of
(\ref{triv42}) .

\subsection{The generalized algebra ${su_{\varphi}}(2)$}

Let $\tt{ F}_{i}, i=0,1$ be the Fock spaces  equipped respectively with  bases
 $\left\{ {\varphi_{n}(x_{i})}\right\}_{n=0}^{\infty }$ and $a^{\pm}_{{\tt F}_{i}},N_{{\tt F}_{i}}, i=0,1$ be the generators
of the generalized oscillators algebra $A_{\varphi }$. These generators
$a^{\pm}_{{\tt F}_{i}},N_{{\tt F}_{i}}, i=0,1$ are  generators of an algebra of the system of the two independent
oscillators if they satisfy the following commutation relations:
\begin{align}
a^{-}_{{\tt F}_{i}}a^{+}_{{\tt F}_{i}}&=2B(N_{{\tt F}_{i}}+I_{{\tt F}_{i}}),&
a^{+}_{{\tt F}_{i}}a^{-}_{{\tt F}_{i}}&=2B(N_{{\tt F}_{i}}),&
[N_{{\tt F}_{i}},a^{\pm}_{{\tt F}_{i}}]&=\pm{a^{\pm}_{{\tt F}_{i}}},\notag\\
[a^{\pm}_{{\tt F}_{1}}&,a^{\pm}_{{\tt F}_{2}}]=0,&
[N_{{\tt F}_{1}},a^{\pm}_{{\tt F}_{2}}]&=0,&
[N_{{\tt F}_{2}},a^{\pm}_{{\tt F}_{1}}]&=0.
\end{align}
We denote by ${su_{\varphi}}(2)$ an algebra generated by the generators  $J^{\varphi}_{+},J^{\varphi}_{-},J^{\varphi}_{z}$,
which are connected with the generators $a^{\pm}_{{\tt F}_{i}},N_{{\tt F}_{i}}$ according to the rules:
\begin{equation}
J^{\varphi}_{+}=a^{+}_{{\tt F}_{1}}a^{-}_{{\tt F}_{2}},\quad
J^{\varphi}_{-}=a^{+}_{{\tt F}_{2}}a^{-}_{{\tt F}_{1}},\quad
J^{\varphi}_{z}=2^{-1}(N_{{\tt F}_{1}}-N_{{\tt F}_{2}}).
\label{triv46}
\end{equation}
\begin{thm}
Let the function $B(x)$ defined by  (\ref{triv40}) is a solution to the following equation:
\begin{equation}
f(x)f(y+1)-f(y)f(x+1)=f(x-y).
\label{triv48}
\end{equation}
Then the operators $J^{\varphi}_{+},J^{\varphi}_{-},J^{\varphi}_{z}$ in the space ${\tt F}_{1}\otimes {\tt F}_{2}$
obey to the following commutation relations:
\begin{equation}
[J^{\varphi}_{z},J^{\varphi}_{\pm}]={{\pm}J^{\varphi}_{\pm}},\qquad
[J^{\varphi}_{+},J^{\varphi}_{-}]=2B(J^{\varphi}_{z}).
\label{triv49}
\end{equation}
\end{thm}
The proof is by direct calculation.
\begin{rem}
We see at once that the relations (\ref{triv49}) are the extensions of the usual commutation relations of
the algebra $su(2)$ and reduce to  the latter in the case $B(x)=x$. An algebra generated by the generators
$J^{\varphi}_{\pm},J^{\varphi}_{z}$ complying with (\ref{triv49}) is called a deformed algebra ${SU}_{q}(2)$
corresponding to the orthonormal system $\left\{ {\varphi_{n}(x)}\right\}_{n=0}^{\infty }$. Indeed, it follows
from the next lemma that  all solutions to the equation (\ref{triv48}) make up an one-parameter family with
the parameter $q$.
\end{rem}
\begin{lem}
If a function $f(x)$ is analytical in the region $|x|<R$, where $R>1$,  and satisfies to the equation (\ref{triv48}),
then one can represent it in the following form:
\begin{equation}
f(x)=\frac{sinh(\eta x)}{sinh(\eta )},\qquad
exp(\eta)=q.
\label{triv50}
\end{equation}
\end{lem}
The proof is left to the reader.
\begin{rem}
1. The following functions:
\begin{equation}
B(x)=x,\qquad
B(x)=[x,q]=\frac{{q^x}-{q^{-x}}}{q-q^{-1}}.
\label{triv51}
\end{equation}
give us some examples of solutions to the equation (\ref{triv48}). Note that  the solution $B(x)=x$ corresponds to
the usual  harmonic oscillator and to the algebra $su(2)$; in the case $B(x)=[x,q]$ we obtain the deformed
oscillator ([8],[9]) and the quantum group ${SU}_{q}(2)$.

2. If $B(x)$ does not a solution to (\ref{triv48}), then the second of commutation relation  (\ref{triv49}) takes the
form
\begin{equation}
[J^{\varphi}_{+},J^{\varphi}_{-}]=(J^{\varphi}_{z})F({J^{\varphi}_{z}},{C_z}),
\label{triv52}
\end{equation}
Here $C_z={N_{1}}+{N_{2}}$ is a element of the center of the algebra ${su_{\varphi}}(2)$ generated
by $J^{\varphi}_{\pm},J^{\varphi}_{z}$ complying with (\ref{triv52}) and the first of relation  (\ref{triv49}).
The function $F$ in the right-hand side (\ref{triv52}) is an analytical function in  its own arguments.
\end{rem}

\subsection{}

In this subsection we will  provide the following answer. How is  a measure $\mu$ to be so that the momentum
operator $P_{{\tt G}_{1}}$ satisfies
\begin{equation}
P_{{\tt G}_{1}}(uv)={vP_{{\tt G}_{1}}(u)}+{uP_{{\tt G}_{1}}(v)}.
\label{triv53}
\end{equation}
Let ${\tt G}_{1}\subset{{\tt H}_{1}}$ and the set  ${{\tt G}_{1}}$ be dense in the space ${{\tt H}_{1}}$.
Denote by $\overline{P}_{{\tt G}_{1}}$ the closure of the momentum operator $P_{{\tt G}_{1}}$
in ${{\tt H}_{1}}$. It is easy to prove the next theorem.
\begin{thm}
\label{su}
The operator $\overline{P}_{{\tt G}_{1}}$ in ${{\tt H}_{1}}$ satisfy (\ref{triv53}) if and only if the following
conditions are held:
\begin{align}
\overline{P}_{{\tt G}_{1}}&={\imath}{\sqrt{2} }{{a^{-}}_{{\tt H}_{1}}},\notag\\
{b_{n}}^{2}&=(n+1){{b_{0}}^{2}},\quad
n\geq 1.
\label{dd}
\end{align}
\end{thm}
\begin{rem}
The second of condition  (\ref{dd}) means that an appropriate oscillator is the usual quantum mechanical one.
\end{rem}

Below we consider the examples of  generalized oscillators algebras corresponding to the classical
orthogonal polynomials.

\section{Hermite polynomials}

\subsection{}

First we consider the main example underlying  our construction, namely, the Hermite polynomials
(\cite{sze},\cite{nus},\cite{koe}).

Let ${{\tt G}_{1}}={L^{2}}(R),$\quad  ${{\tt H}_{1}}={L^{2}}(R;{\frac{1}{\sqrt{\pi }}}{\exp(-x^{2})}{dx}$ and
\begin{equation}
\psi (x)= {\pi ^{-\frac{1}{4}}}\exp(-\frac{x^{2}}{2}).
\label{triv56}
\end{equation}
We denote by ${H_{n}}(x)$ the Hermite polynomials
\begin{equation}
{H_{n}}(x)= n!\sum_{\nu =0}^{[\frac{n}{2}]}{\frac{(-1)^{\nu }}{{\nu }!}}{\frac{(2x)^{n-2\nu }}{(n-2\nu )!}},
\label{triv57}
\end{equation}
We define the functions $\left\{ {\psi_{n}(x)}\right\}_{n=0}^{\infty }$ and $\left\{ {\phi_{n}(x)}\right\}_{n=0}^{\infty }$
by the following formulas:
\begin{equation}
{\psi_{n}}(x)= {\sqrt[4]{\pi }}{{d_{n}}^{-1}}{H_{n}}(x),\qquad
{\phi_{n}}(x)={{d_{n}}^{-1}}{\exp(-\frac{x^{2}}{2})}{H_{n}}(x),\qquad
n\geq 0,
\label{triv59}
\end{equation}
where
\begin{equation}
d_{n}= ({2^{n}}{n!}{\sqrt \pi } )^{\frac{1}{2}},\qquad
n\geq 0.
\label{triv58}
\end{equation}
The recurrence relations for the Hermite polynomials (\cite{sze}) give us the formulas (\ref{triv5}) ,(\ref{triv6}) with
\begin{equation}
b_{n}= {\frac{1}{2}}(\frac{d_{n+1}}{d_{n}})=\sqrt{\frac{n+1}{2}} .
\label{triv60}
\end{equation}
 From the Mehler formula for the Hermite polynomials (\cite{meh}) the following expression
for  the Poisson kernel follows:
\begin{equation}
\label{triv61}
{{\pi} ^{-\frac{1}{2}}}{\sum_{n=0}^{\infty }{{\omega ^{n}}\cdot {\psi_{n}(x)}\cdot {\psi_{n}(y)}}}=
{(1-{\omega }^{2})^{-\frac{1}{2}}}
\exp(\frac{{2xy\omega }-({x^{2}}+{y^{2}}){{\omega }^{2}}}{1-{\omega }^{2}}).
\end{equation}
Combining (\ref{triv61}) with the definition of the (direct and inverse) generalized Fourier transform conforming to
the orthonormal system $\left\{ {\varphi_{n}(x)}\right\}_{n=0}^{\infty }$ we get
\begin{equation}
F_{\phi}={K_{\tt G}}(-\imath),\qquad
{F_{\phi}}^{-1}={K_{\tt G}}(\imath),
\end{equation}
where respectively
\begin{equation}
{\mathfrak{K}_{\tt F}(x,y;-\imath)}=\frac{\exp({{-\imath}xy})}{\sqrt {2\pi }},\qquad
{\mathfrak{K}_{\tt F}(x,y;\imath))}=\frac{\exp({\imath}xy)}{\sqrt {2\pi }}.
\label{triv62}
\end{equation}
Let us remark that in this case the generalized  Fourier transform be the same as the usual Fourier transform.
An easy computation shows that we have in the space ${{\tt H}_{1}}$:
\begin{equation}
\overline{P_{{\tt G}_{1}}}={\imath}\frac{d}{dx},
\end{equation}
that is conforming to the theorem  {\ref{su}} since the conditions (\ref{dd}) are valid.
Note also that
\begin{align}
{a^{+}}_{{\tt G}_{1}}&={\frac{1}{\sqrt{2}}}
\left (
{X_{{\tt G}_{1}}}-{\frac{d}{dx}}
\right )
,\qquad
{a^{-}}_{{\tt G}_{1}}={\frac{1}{\sqrt{2}}}
\left (
{X_{{\tt G}_{1}}}+{\frac{d}{dx}}
\right ),\notag\\
H_{{\tt G}_{1}}&=(X_{{\tt G}_{1}})^{2}+({P_{{\tt G}_{1}}})^{2}=
{-\frac{d^2}{{dx}^2}}+x^2.
\end{align}
Then the equation
\begin{equation}
{H_{{\tt G}_{1}}}{{\phi_{n}}(x)}={\lambda _{n}}{\phi_{n}}(x),\quad
\lambda _{n}=2n+1,
\end{equation}
takes the following form:
\begin{equation}
({-\frac{d^2}{{dx}^2}}+x^2){{\phi_{n}}(x)}=(2n+1){{\phi_{n}}(x)}.
\label{triv63}
\end{equation}
It can easily be checked that (\ref{triv63}) is equivalent to the well-known equation for the Hermite polynomials:
\begin{equation}
{{{H_{n}}^{\prime \prime }}(x)}-{2x{{H_{n}}^{\prime}}(x)}+{2n{H_{n}}(x)}=0.
\end{equation}
In the next section we present  the first substantive example, namely, the ultraspherical polynomials.

\section{ Ultrasferical polynomials }

First we consider a particular case of the ultraspherical polynomials, namely, the Legendre polynomials.

\subsection{The Legendre polynomials}

Let
$$
{{\tt G}_{1}}={L^{2}}([-1,1]),\quad
{{\tt H}_{1}}={L^{2}}([-1,1];2^{-1}),
$$
and the function  $\psi (x)=\frac{1}{\sqrt{2}}$.
The Legendre polynomials are defined by
\begin{equation}
{P_{n}}(x)={}_{2}F_{1}
\left (
-n,n+1;1;\frac{1-x}{2}
\right ).
\label{triv64}
\end{equation}
The functions of the orthonormal systems  $\left\{ {\psi_{n}(x)}\right\}_{n=0}^{\infty }$
and  $\left\{ {\phi_{n}(x)}\right\}_{n=0}^{\infty }$
are given by the following formulas:
\begin{equation}
{\psi_{n}}(x)= {\sqrt {2}}{\phi_{n}}(x),\qquad
{\phi_{n}}(x)=\sqrt{\frac{2n+1}{2}} {P_{n}}(x),\qquad
n\geq 0.
\label{triv65}
\end{equation}
Taking into account the recurrence relations for the Legendre polynomials (\cite{sze}) we obtain the formulas
(\ref{triv5}) ,(\ref{triv6}) where
\begin{equation}
b_{n}=\sqrt{\frac{(n+1)^{2}}{(2n+1)(2n+3)}},\quad
n\geq 0,
\label{triv66}
\end{equation}
and
\begin{equation}
{b_{n}^{2}}-{b_{n-1}^{2}}= -{(2n-1)^{-1}}{(2n+1)^{-1}}{(2n+3)^{-1}} .
\label{triv67}
\end{equation}
In the construction of the momentum operator we will use the following  differential operator:
\begin{equation}
A=(1-x^{2})\frac{d}{dx},
\label{triv68}
\end{equation}
in the space ${\tt G}_{1}$. The operator $A$ acts on the basis vectors $\left\{ {\phi_{n}(x)}\right\}_{n=0}^{\infty }$
of  ${\tt G}_{1}$ by:
\begin{equation}
A{\phi_{n}}(x)= {(n+1){b_{n-1}}{\phi_{n-1}}(x)}-
{n{b_{n}}{\phi_{n+1}}(x)},\qquad
n\geq 0,\quad
b_{-1}=0.
\label{triv69}
\end{equation}
Using the definition (\ref{triv38}) of the number operator $N$ and  (\ref{triv35}), from (\ref{triv5}) ,(\ref{triv6})
and  (\ref{triv68}) ,(\ref{triv69})  we  have
\begin{eqnarray}
P_{{\tt H}_{1}}={ \imath}{(N_{{\tt H}_{1}}-{(2^{-1}){I_{{\tt H}_{1}}}})^{-1}}
{(N_{{\tt H}_{1}}+3{(2^{-1}){I_{{\tt H}_{1}}}})^{-1}}\nonumber\\
({(N_{{\tt H}_{1}}+{(2^{-1}){I_{{\tt H}_{1}}}})A}
-{{X_{{\tt H}_{1}}}{N_{{\tt H}_{1}}}}-
(2^{-1})((N_{{\tt H}_{1}}+{3(2^{-1})}{I_{{\tt H}_{1}})}{X_{{\tt H}_{1}}}).
\label{triv70}
\end{eqnarray}
It can easily be checked that the formula (\ref{triv29}) for the operator $P_{{\tt H}_{1}}$ is valid:
\begin{equation}
{P_{{\tt H}_{1}}}{\phi_{n}}(x)=\imath({{b_{n-1}}{\phi_{n-1}}(x)}-
{{b_{n}}{\phi_{n+1}}(x)}),\qquad
n\geq 0,\quad
b_{-1}=0.
\label{triv71}
\end{equation}
Then the ladder operators $a^{-}_{{\tt H}_{1}}$ and  $a^{+}_{{\tt H}_{1}}$ are given by:
\begin{eqnarray}
a^-_{{\tt H}_{1}}=\frac{1}{\sqrt{2}}(N_{{\tt H}_{1}}+3(2^{-1})I_{{\tt H}_{1}})^{-1}(A+X_{{\tt H}_{1}}N_{{\tt H}_{1}}),\nonumber\\
a^+_{{\tt H}_{1}}=\frac{1}{\sqrt{2}}(N_{{\tt H}_{1}}-(2^{-1})I_{{\tt H}_{1}})^{-1}
(-A+X_{{\tt H}_{1}}(N_{{\tt H}_{1}}+I_{{\tt H}_{1}})),.
\label{triv72}
\end{eqnarray}
Further, the eigenvalue of the operator $H_{{\tt H}_{1}}=(X_{{\tt H}_{1}})^{2}+({P_{{\tt H}_{1}}})^{2}$
(the energy levels) amount to:
\begin{equation}
{\lambda _{0}}={\frac{2}{3}},\quad
{\lambda _{n}}=
{\frac{n(n+1)-{2^{-1}}}{(n+3(2^{-1}))(n-(2^{-1}))}},\quad
n>0.
\label{triv73}
\end{equation}
\begin{thm}
\label{dod}
The equation $H_{{\tt H}_{1}}{\phi_{n}}(x)={\lambda _{n}}{\phi_{n}}(x),\quad n\geq 0$, where
 the eigenvalue  $\lambda _{n}$ of the operator $H_{{\tt H}_{1}}=(X_{{\tt H}_{1}})^{2}+({P_{{\tt H}_{1}}})^{2}$
defined by (\ref{triv73}),  is equivalent to the usual differential equation for  the Legendre polynomials ([12]) :
\begin{equation}
{\frac{d}{dx}((1-x^{2})\frac{d}{dx}){{P_{n}}(x)}}+{n(n+1){{P_{n}}(x)}}=0,\quad
n\geq 0.
\label{triv74}
\end{equation}
\end{thm}
\begin{proof}
On  account  of (\ref{triv73}) we rewrite the equation
$H_{{\tt H}_{1}}{\phi_{n}}(x)={\lambda _{n}}{\phi_{n}}(x),\quad n\geq 0$, as an operator equality in the
space ${\tt H}_{1}$:
\begin{eqnarray}
(X_{{\tt H}_1})^2+(P_{{\tt H}_1})^2=(N_{{\tt H}_1}(N_{{\tt H}_1}+I_{{\tt H}_1})-(2^{-1})I_{{\tt H}_{1}})\nonumber\\
(N_{{\tt H}_1}-(2^{-1})I_{{\tt H}_1})^{-1}(N_{{\tt H}_1}+3(2^{-1})I_{{\tt H}_{1}})^{-1}.
\label{triv75}
\end{eqnarray}
Substituting (\ref{triv70}) in (\ref{triv75}) we get
\begin{align}
(2^{-1})I_{{\tt H}_{1}}&+{\imath}((N_{{\tt H}_{1}}+(2^{-1})I_{{\tt H}_{1}})A
-X_{{\tt H}_{1}}N_{{\tt H}_{1}}-\notag\\
(2^{-1})((N_{{\tt H}_{1}}&+3(2^{-1})I_{{\tt H}_{1}})X_{{\tt H}_{1}})P_{{\tt H}_{1}}+
(N_{{\tt H}_{1}}-(2^{-1})I_{{\tt H}_{1}})\notag\\
(N_{{\tt H}_{1}}&+3(2^{-1})I_{{\tt H}_{1}})
X_{{\tt H}_{1}}^{2}=N_{{\tt H}_{1}}(N_{{\tt H}_{1}}+I_{{\tt H}_{1}}).
\label{triv76}
\end{align}
It is not hard to prove that
 \begin{align}
-A^{2}+{X_{{\tt H}_{1}}}^{2}N_{{\tt H}_{1}}(N_{{\tt H}_{1}}+I_{{\tt H}_{1}})&=
(2^{-1})I_{{\tt H}_{1}}+{\imath}((N_{{\tt H}_{1}}+(2^{-1})I_{{\tt H}_{1}})A\notag\\
-X_{{\tt H}_{1}}N_{{\tt H}_{1}}-(2^{-1})((N_{{\tt H}_{1}}&+3(2^{-1})I_{{\tt H}_{1}})X_{{\tt H}_{1}})
P_{{\tt H}_{1}}+\notag\\
(N_{{\tt H}_{1}}-(2^{-1})I_{{\tt H}_{1}})(N_{{\tt H}_{1}}&+3(2^{-1})I_{{\tt H}_{1}})
X_{{\tt H}_{1}}^{2}.
\label{triv77}
\end{align}
Then from (\ref{triv70}) and (\ref{triv77}) we have
\begin{equation}
({A^{2}}+
({I_{{\tt H}_{1}}}-{X_{{\tt H}_{1}}^{2}})
n(n+1)){\phi_{n}}(x)=0.
\label{triv79}
\end{equation}
Substituting (\ref{triv68}) and (\ref{triv65}) in (\ref{triv79}) we get the equation (\ref{triv74}).
\end{proof}
\begin{rem}
Using (\ref{triv79}) one can get
\begin{equation}
(N_{{\tt H}_{1}}+{(2^{-1}){I_{{\tt H}_{1}}}})=D=(2^{-1})
({I_{{\tt H}_{1}}}-4{\frac{d}{dx}}
((1-x^{2})\frac{d}{dx}))^{\frac{1}{2}}.
\label{triv80}
\end{equation}
We exclude the number operator $N_{{\tt H}_{1}}$ from the right-side of (\ref{triv70}).
Then we obtain
\begin{align}
\label{triv81}
P_{{\tt H}_{1}}&=
{ \imath}
{({D^{2}}-{I_{{\tt H}_{1}}})^{-1}}
(DA-{X_{{\tt H}_{1}}}D
-{(2^{-1})D
{X_{{\tt H}_{1}}}}),\\
\label{triv82}
{a^{-}}_{{\tt H}_{1}}&={\frac{1}{\sqrt{2} }}
{(D+{I_{{\tt H}_{1}}})^{-1}}
(A+{{X_{{\tt H}_{1}}}(D-(2^{-1}){I_{{\tt H}_{1}}})}),\\
\label{triv83}
{a^{+}}_{{\tt H}_{1}}&={\frac{1}{\sqrt{2} }}
{(D-{I_{{\tt H}_{1}}})^{-1}}
({-A}+{{X_{{\tt H}_{1}}}(D+(2^{-1}){I_{{\tt H}_{1}}})}).
\end{align}
\end{rem}
Furthermore
\begin{align}
\label{triv84}
H_{{\tt H}_{1}}&=
({D^{2}}-{\frac{3}{4}}{I_{{\tt H}_{1}}})
{({D^{2}}-{I_{{\tt H}_{1}}})^{-1}}=\notag\\
{I_{{\tt H}_{1}}}+{\frac{1}{4}}{({D^{2}}-{I_{{\tt H}_{1}}})^{-1}}&=
{I_{{\tt H}_{1}}}-(3{I_{{\tt H}_{1}}}+{\frac{d}{dx}}
((1-x^{2})\frac{d}{dx}))^{-1}.
\end{align}
and in view of (\ref{triv67}) the commutation relations (\ref{triv42}) for the operators  (\ref{triv82})-(\ref{triv83})
will look like:
\begin{equation}
[{a^{-}}_{{\tt H}_{1}},{a^{+}}_{{\tt H}_{1}}]={-\frac{1}{4}}
{({{N_{{\tt H}_{1}}}^{2}}-{\frac{1}{4}}{I_{{\tt H}_{1}}})^{-1}}
{({{N_{{\tt H}_{1}}}^{2}}+{\frac{3}{2}}{I_{{\tt H}_{1}}})^{-1}}.
\label{triv85}
\end{equation}

Now we turn to the general case of the ultraspherical polynomials.

\subsection{The Gegenbauer polynomials}
Let
$$
{{\tt G}_{1}}={L^{2}}([-1,1]),
$$
$$
{{\tt H}_{1}}={L^{2}}([-1,1];{({d_{0}}(\alpha ))^{-2}}{(1-{x^{2}})^{\alpha }}dx),
$$
where
$$
{{d_{0}}^{2}}(\alpha )={2^{{2\alpha }+1}}\frac{(\Gamma (\alpha +1))^{2}}{\Gamma (2(\alpha +1))}.
$$
The ultraspherical polynomials are defined by  the hypergeometric function (\cite{bet, pbm}):
\begin{equation}
{P_{n}^{(\alpha ,\alpha )}}(x)={\frac{(\alpha +1)_{n}}{n!}}
{{}_{2}F_{1}
\left (
-n,n+2\alpha +1;\alpha +1;\frac{1-x}{2}
\right )}.
\label{triv86}
\end{equation}
The Pochhammer-symbol (\cite{gas}) $(\beta )_{n}$ is defined by ${(\beta )_{0}}=1,\quad {(\beta )_{n}}=\beta (\beta +1)\cdots(\beta +n-1),
\quad n\geq 1$.
For  $\alpha >-1$ the following orthogonal relations are valid:
$$
\int_{-1}^{1}{{{P_{n}^{(\alpha ,\alpha )}}(x)}
{{P_{m}^{(\alpha ,\alpha )}}(x)}{(1-x^{2})^{\alpha }}{dx}}=
{{d_{n}}^{2}}{\delta _{mn}},\quad
n,m\geq 0,
$$
with the constant of normalization $d_{n}$ given by
\begin{equation}
{d_{n}}^{2}=\frac{{2^{2\alpha +1}}{(\Gamma (n+\alpha +1))^{2}}}{(2n+2\alpha +1){n!}\Gamma (n+2\alpha +1)},\quad
n\geq 0.
\label{triv88}
\end{equation}
The Gegenbauer polynomials are defined as usual  (\cite{sze}):
\begin{equation}
{P_{n}^{(\lambda )}}(x)
={\frac{{\Gamma (\alpha +1)}{\Gamma (n+{2\alpha} +1)}}
{{\Gamma (2\alpha +1)}{\Gamma (n+\alpha +1)}}}
{P_{n}^{(\alpha ,\alpha)}}(x),\qquad
\alpha =\lambda -{2^{-1}},\qquad
\label{triv90}
\end{equation}
($\lambda >-{2^{-1}},\quad
n\geq 0$).

Let $\psi (x)={{d_{0}}^{-1}}(1-x^{2})^{{2^{-1}}\alpha }$. We determine a functions of the orthonormal systems
$\left\{ {\psi_{n}(x)}\right\}_{n=0}^{\infty }$  and  $\left\{ {\phi_{n}(x)}\right\}_{n=0}^{\infty }$  by the following formulas:
\begin{equation}
{\psi_{n}}(x)={d_{0}}{{d_{n}}^{-1}}{P_{n}^{(\alpha ,\alpha)}}(x),\qquad
{\phi_{n}}(x)={\psi(x)}{\psi_{n}}(x),\qquad
n\geq 0,
\label{triv89}
\end{equation}
where $d_{n}$ is given via (\ref{triv88}).
The function ${\psi_{n}}(x)$ defined by (\ref{triv86}) satisfies the relations (\ref{triv5}) and (\ref{triv6}), where
\begin{equation}
b_{n}=\sqrt{\frac{(n+1)(n+2\alpha +1)}{(2n+2\alpha +1)(2n+2\alpha +3)}},\quad
n\geq 0,\quad
b_{-1}=0,
\label{triv91}
\end{equation}
and
\begin{equation}
{{b_{n}}^{2}}-{{b_{n-1}}^{2}}= {2^{-1}}\lambda (\lambda -1)
{(n+\lambda )^{-1}}{(n-1+\lambda )^{-1}}{(n+\lambda +1)^{-1}} .
\label{triv92}
\end{equation}
In order to find a differential expression for the momentum operator $P_{{\tt H}_{1}}$
we use the known formula (\cite{sze}):
\begin{equation}
A{\psi_{n}}(x)= {(n+2\alpha +1){b_{n-1}}{\psi_{n-1}}(x)}-
{n{b_{n}}{\psi_{n+1}}(x)},\qquad
n\geq 0,\quad
b_{-1}=0,
\label{triv93}
\end{equation}
where $A$ and $b_{n}$ are defined by (\ref{triv68}) and (\ref{triv89}) respectively.
Combining (\ref{triv91}) and (\ref{triv5}) with (\ref{triv89}) we get
\begin{align}
a^{-}_{{\tt H}_{1}}&={\frac{1}{\sqrt{2} }}(N_{{\tt H}_{1}}+(\alpha +3{(2^{-1})){I_{{\tt H}_{1}}}})^{-1}
(A+X_{{\tt H}_{1}}N_{{\tt H}_{1}}),\notag\\
a^{+}_{{\tt H}_{1}}&={\frac{1}{\sqrt{2} }}(N_{{\tt H}_{1}}+(\alpha -{(2^{-1})){I_{{\tt H}_{1}}}})^{-1}
(-A+X_{{\tt H}_{1}}(N_{{\tt H}_{1}}+(2\alpha +1)I_{{\tt H}_{1}})),
\label{triv94}
\end{align}
which generalize (\ref{triv72}) and reduce to these as $\alpha =0$.
From (\ref{triv94}) and the formula $P_{{\tt H}_{1}}={\frac{1}{{\imath}\sqrt{2} }}({a^{+}}_{{\tt H}_{1}}
-{a^{-}}_{{\tt H}_{1}})$ it follows that:
\begin{eqnarray}
P_{{\tt H}_{1}}={ \imath}{(N_{{\tt H}_{1}}+(\alpha -{(2^{-1})){I_{{\tt H}_{1}}}})^{-1}}
({(N_{{\tt H}_{1}}+(\alpha +3{(2^{-1})){I_{{\tt H}_{1}}}})^{-1}}\nonumber\\
{(N_{{\tt H}_{1}}+(\alpha +{(2^{-1})){I_{{\tt H}_{1}}}})A}
-{(N_{{\tt H}_{1}}+(\alpha +3{(2^{-1})){I_{{\tt H}_{1}}}})^{-1}}\nonumber\\
{{X_{{\tt H}_{1}}}{N_{{\tt H}_{1}}}}-
(\alpha +(2^{-1})){X_{{\tt H}_{1}}}).
\label{triv95}
\end{eqnarray}
Note that (\ref{triv95}) generalize (\ref{triv70}) for the Legendre polynomials and  reduce to these as $\alpha =0$.
\begin{rem}
The energy operator $H_{{\tt H}_{1}}=(X_{{\tt H}_{1}})^{2}+({P_{{\tt H}_{1}}})^{2}$ is bounded and has the energy
 levels
\begin{equation}
{\lambda _{0}}={\frac{2}{2\alpha +3}},\quad
{\lambda _{n}}=
{\frac{n(n+2\alpha +1)+(\alpha -{2^{-1}})}{(n+\alpha +3(2^{-1}))(n+\alpha -(2^{-1}))}},\quad
n>0.
\label{triv96}
\end{equation}
\end{rem}
The next theorem is the extension of the analogous theorem  {\ref{dod}}.
\begin{thm}
\label{dof}
The equation $H_{{\tt H}_{1}}{\phi_{n}}(x)={\lambda _{n}}{\phi_{n}}(x),\quad n\geq 0$, where
 the eigenvalues  $\lambda _{n}$ of the operator $H_{{\tt H}_{1}}=(X_{{\tt H}_{1}})^{2}+({P_{{\tt H}_{1}}})^{2}$  are
defined by (\ref{triv96}),  is equivalent to the usual differential equation for  the ultraspherical  polynomials ([12]) :
\begin{eqnarray}
(1-x^2)\frac{d}{dx}((1-x^2)^{\alpha -1}\frac{d}{dx})P_n^{(\alpha ,\alpha )}(x)+\nonumber\\
n(n+2\alpha +1)(1-x^2)^{\alpha +1}P_n^{(\alpha ,\alpha )}(x)=0,
\label{triv97}
\end{eqnarray}
$(n\geq 0)$.
\end{thm}

A slight change in the proof of  the theorem {\ref{dod}} shows that the theorem  {\ref{dof}} is true.
\begin{rem}
Similarly to (\ref{triv81})-(\ref{triv83}) we have the following formulas:
\begin{align}
P_{{\tt H}_{1}}&=
{ \imath}
{({D_{{\tt H}_{1}}}^{2}-{I_{{\tt H}_{1}}})^{-1}}
(D_{{\tt H}_{1}}A-{X_{{\tt H}_{1}}}D_{{\tt H}_{1}}
-{(\alpha +2^{-1})D_{{\tt H}_{1}}
{X_{{\tt H}_{1}}}}),\\
a^{-}_{{\tt H}_{1}}&={\frac{1}{\sqrt{2} }}
{(D_{{\tt H}_{1}}+{I_{{\tt H}_{1}}})^{-1}}
(A+{{X_{{\tt H}_{1}}}(D_{{\tt H}_{1}}-(\alpha +2^{-1}){I_{{\tt H}_{1}}})}),\\
a^{+}_{{\tt H}_{1}}&={\frac{1}{\sqrt{2} }}
{(D_{{\tt H}_{1}}-{I_{{\tt H}_{1}}})^{-1}}
({-A}+{{X_{{\tt H}_{1}}}(D_{{\tt H}_{1}}+(\alpha +2^{-1}){I_{{\tt H}_{1}}})}).
\end{align}
 Analogy with (\ref{triv84}) gives us:
\begin{align}
H_{{\tt H}_{1}}&=
({{D_{{\tt H}_{1}}}^{2}}-({\alpha }^{2}+\frac{3}{4}){I_{{\tt H}_{1}}})
{({{D_{{\tt H}_{1}}}^{2}}-{I_{{\tt H}_{1}}})^{-1}}\notag\\
{}&={I_{{\tt H}_{1}}}+(\frac{1}{4}-{\alpha }^{2}){({{D_{{\tt H}_{1}}}^{2}}-{I_{{\tt H}_{1}}})^{-1}}\notag\\
{}&={I_{{\tt H}_{1}}}+(\frac{1}{4}-{\alpha }^{2})({({-\frac{3}{4}}+{{\alpha }^{2}}(1-x^{2})^{-1})I_{{\tt H}_{1}}}-{\frac{d}{dx}}
((1-x^{2})\frac{d}{dx}))^{-1}.
\end{align}
In view of (\ref{triv92}) the commutation relations (\ref{triv42}) for the operators  (\ref{triv94})
will look like:
\begin{equation}
[a^{-}_{{\tt H}_{1}},a^{+}_{{\tt H}_{1}}]=\lambda (\lambda -1)
{(({N_{{\tt H}_{1}}+\alpha {I_{{\tt H}_{1}}})^{2}}-{\frac{1}{4}}{I_{{\tt H}_{1}}})^{-1}}
{({N_{{\tt H}_{1}}}+(\alpha +\frac{3}{2}){I_{{\tt H}_{1}}})^{-1}}.
\end{equation}
\end{rem}
In conclusion of this section we consider yet another special case, namely, the Chebyshev polynomials.

\subsection{Chebyshev polynomials}

The Chebyshev polynomials of the first kind ${T_{n}}(x)$ and those of the second kind ${U_{n}}(x)$
are special cases of the Gegenbauer polynomials for $\lambda =0$($\alpha =-{2^{-1}}$) and
$\lambda =1$ ($\alpha =2^{-1}$) respectively. In both cases,  it follows from (\ref{triv92}):
$$
b_{-1}=0,\quad
b_{n}=\frac{1}{2},\quad
n\geq 0.
$$
Hence  all operators - co-ordinate, momentum and hamiltonian - are the same  for
the Chebyshev polynomials of the first kind, as for those of the second kind. Thus both of
these polynomials systems give us the unitary equivalent representations of  the same
oscillator in the different spaces:
 $$
{\tt H}_{1}=L^{2}([-1,1];(\sqrt{1-x^{2}})^{-1}\frac{dx}{{{d_{0}}^{2}}(2^{-1})}),
$$
and
$$
{\tt H}_{2} =L^{2}([-1,1];{\sqrt{1-x^{2}}}\frac{dx}{{{d_{0}}^{2}}(-2^{-1})}).
$$
In both cases the energy levels are equal:
$$
\lambda _{0}=\frac{1}{2},\quad
\lambda _{n}=1,\quad
n\geq 1.
$$
Now we consider the following orthonormal systems in the spaces ${\tt H}_{1}$ and ${\tt H}_{2}$
respectively:
\begin{align}
\psi_{n}^{(1)}&={\sqrt{2n}}{\frac{{\Gamma (n+\frac{1}{2})}^{2}}
{{n!}\Gamma (n)}}{2^{-2n}}{C_{2n}^{n}}{T_{n}}(x),\\
{T_{n}}(x)&={2^{2n}}{{C_{2n}^{n}}^{-1}}{{P_{n}}^{(-2^{-1},-2^{-1})}}(x),\\
\psi_{n}^{(2)}&={\sqrt{2n+2}}{\frac{{\Gamma (n+3\frac{1}{2})}^{2}}
{{n!}\Gamma (n+2)}}{2^{-2(n+1)}}{C_{2(n+1)}^{n+1}}{U_{n}}(x),\\
{U_{n}}(x)&={2^{2n+1}}{{C_{2(n+1)}^{n+1}}^{-1}}{{P_{n}}^{(2^{-1},2^{-1})}}(x),
\end{align}
Then the operator $A$ acts on basis vectors in the mentioned spaces by
\begin{align}
A{\psi_{n}^{(1)}}&= n{b_{n-1}}{\psi_{n-1}^{(1)}}-
n{b_{n}}{\psi_{n+1}^{(1)}},\qquad
n\geq 0,\quad
b_{-1}=0.\\
A{\psi_{n}^{(2)}}&= (n+2){b_{n-1}}{\psi_{n-1}^{(2)}}-
n{b_{n}}{\psi_{n+1}^{(2)}},\qquad
n\geq 0,\quad
b_{-1}=0.
\end{align}
As before, we have
\begin{align}
a^{-}_{{\tt H}_{1}}&={\frac{1}{\sqrt{2} }}
{({N_{{\tt H}_{1}}}+{I_{{\tt H}_{1}}})^{-1}}
(A+{{X_{{\tt H}_{1}}}{N_{{\tt H}_{1}}}}),\\
a^{+}_{{\tt H}_{1}}&={\frac{1}{\sqrt{2} }}
{(N_{{\tt H}_{1}}-{{I_{{\tt H}_{1}}}})^{-1}}
(-A+{{X_{{\tt H}_{1}}}{N_{{\tt H}_{1}}}}),\\
a^{-}_{{\tt H}_{2}}&={\frac{1}{\sqrt{2} }}
{({N_{{\tt H}_{2}}}+{I_{{\tt H}_{2}}})^{-1}}
(A+{{X_{{\tt H}_{2}}}{N_{{\tt H}_{2}}}}),\\
a^{+}_{{\tt H}_{2}}&={\frac{1}{\sqrt{2} }}
{N_{{\tt H}_{2}}}^{-1}
(-A+{X_{{\tt H}_{2}}}({N_{{\tt H}_{2}}}+2{I_{{\tt H}_{2}}}),
\end{align}
One can present the momentum operators in the following forms:
\begin{eqnarray}
P_{{\tt H}_1}=
\imath
(N_{{\tt H}_1}-I_{{\tt H}_1})^{-1}
((N_{{\tt H}_1}+I_{{\tt H}_1})^{-1}N_{{\tt H}_1}A
-(N_{{\tt H}_1}+I_{{\tt H}_1})^{-1}
X_{{\tt H}_1}N_{{\tt H}_1}). \\
P_{{\tt H}_2}=\imath (N_{{\tt H}_2})^{-1}
((N_{{\tt H}_2}+
2I_{{\tt H}_2})^{-1}
(N_{{\tt H}_2}+I_{{\tt H}_2})A
-(N_{{\tt H}_2}+2I_{{\tt H}_2})^{-1}
X_{{\tt H}_2}N_{{\tt H}_2}-X_{{\tt H}_2}).
\label{triv101}
\end{eqnarray}

\section{Nonsymmetric Jacobi  matrix}

 Let $\mu$ be a symmetric probability measure on  $R^1$. We construct in this section a noncanonical
orthogonal polynomial system
$\left\{ {\varphi_{n}(x)}\right\}_{n=0}^{\infty }$ in the space ${\tt H}_{1}$. As before, we define $\mu_{2k}$
by (\ref{triv2}) and look for the positive sequences $\left\{ {b_{n}}\right\}_{n=0}^{\infty }$,
$\left\{ {c_{n}}\right\}_{n=0}^{\infty }$ as  solutions of the following equation system:
\begin{equation}
\mu_{0}=1,\quad
\mu_{2k}={{b _{0}}{c _{1}}\cdot( {b _{0}}{c_{1}}+{b _{1}}{c_{2}})
\cdots \sum_{j=0}^{k-1}{{b _{j}}{c_{j+1}}}},\qquad
k= 0,1,\dots .
\label{triv102}
\end{equation}
Contrary to (\ref{triv3}) there is an infinite number of solutions of the system (\ref{triv102}). We can find
uniquely from (\ref{triv102}) only
\begin{equation}
d_{j}=\sqrt{{b _{j}}{c_{j+1}}},\qquad
j= 0,1,\dots .
\label{triv103}
\end{equation}
Now we determine a polynomials system $\left\{ {\varphi_{n}(x)}\right\}_{n=0}^{\infty }$ from the given
sequences $\left\{ {b_{n}}\right\}_{n=0}^{\infty }$, $\left\{ {c_{n}}\right\}_{n=0}^{\infty }$ by the following
recurrence relations:
\begin{equation}
{x {\varphi_{n}(x)}}= {b_{n}{\varphi_{n+1}(x)}}+{c_{n}{\varphi_{n-1}(x)}},\qquad
n\geq 0,\quad
b_{-1}=0,
\label{triv104}
\end{equation}

\begin{equation}
{\varphi_{0}(x)}=1.
\label{triv105}
\end{equation}
The canonical  polynomial system $\left\{ {\psi_{n}(x)}\right\}_{n=0}^{\infty }$ is determined by the
recurrence relations (\ref{triv5}), (\ref{triv6}) replacing $b_{n}$ with $d_{n}$.
From the theorem \ref{a} it follows that the set $\left\{ {\psi_{n}(x)}\right\}_{n=0}^{\infty }$ is a
orthonormal system in the space ${\tt H}_{1}$.
It is clear that  making the following renormalization of the system $\left\{ {\varphi_{n}(x)}\right\}_{n=0}^{\infty }$
\begin{equation}
{\varphi_{n}}(x)={\gamma _{n}}{\psi_{n}}(x),\qquad
n\geq 0,
\label{triv106}
\end{equation}
with
\begin{equation}
{\gamma _{0}}=1,\quad
{\gamma _{n}}=\sqrt{\frac{c_{1}\cdot c_{2}\cdots c_{n} }{b_{0}\cdot b_{2}\cdots b_{n-1}}},\quad
n\geq 1,
\label{107}
\end{equation}
 we get  from (\ref{triv104})  the relations (\ref{triv5}), (\ref{triv6}) replacing $b_{n}$ with $d_{n}$.
Taking into account the orthonormal condition  for  $\left\{ {\psi_{n}(x)}\right\}_{n=0}^{\infty }$ we
obtain the following orthogonal condition  for  $\left\{ {\varphi_{n}(x)}\right\}_{n=0}^{\infty }$:
\begin{equation}
\int_{-\infty }^{\infty }{{\varphi_{n}(x)}{{\varphi_{m}(x)}}{\mu (dx)}}={\gamma _{n}}{\gamma _{m}}
\int_{-\infty }^{\infty }{{\psi_{n}(x)}{{\psi_{m}(x)}}{\mu (dx)}}=
={\gamma _{n}}^{2}{\delta _{nm}}.
\label{108}
\end{equation}
Therefore the set $\left\{ {\psi_{n}(x)}\right\}_{n=0}^{\infty }$
is an orthogonal system, however it does not need to be a orthonormal system
(as $\gamma _{n}^{2}\not= 1$).
\begin{rem}
\label{gg}
This argument shows that for a symmetric probability measure one can reduce the recurrence relations
(\ref{triv104})  to the symmetrical ones (\ref{triv5}).
\end{rem}

\section{Nonsymmetric scheme}

\subsection{}

Let $\mu$ be a probability but not necessarily  a symmetric  measure on $[a,b]\subset R$, i.e. the conditions
(\ref{triv1}) are incorrect. Denote by \begin{math}{\tt H}=L^2([a,b];{\mu)}\end{math} the Hilbert space of the
square-integrable functions with respect to the measure $\mu$ on $[a,b]\subset R$.
Let
\begin{equation}
\mu_{0}=1,\qquad
\mu_{k}=\int_{-\infty}^{\infty}{x^{k}\mu(dx)},\qquad
k= 0,1,\dots ,
\label{trivs1}
\end{equation}
then look for the real sequences $\left\{ {b_{n}}\right\}_{n=0}^{\infty }$,
$\left\{ {a_{n}}\right\}_{n=0}^{\infty }$  as  solutions of the following equation system:
\begin{equation}
A_{k,n}= {b_{n}{A_{k-1,n+1}}}+{a_{n}{A_{k-1,n}}}+{b_{n-1}{A_{k-1,n-1}}},\qquad
n\geq 0,\qquad
b_{-1}=0,
\label{trivs2}
\end{equation}
also satisfying the conditions:
\begin{equation}
A_{0,0}=1,\quad
A_{k,0}=\mu_{k},\quad
A_{0,k}=0,\qquad
k\geq 1.
\label{trivs3}
\end{equation}
\begin{lem}
There is a unique solution to the system of equations (\ref{trivs2}), (\ref{trivs3}) with respect to
the variables $(a_{n},b_{n}A_{k,n}),n\geq 0,k\geq 0$.
\end{lem}
If  sequences $\left\{ {a_{n}}\right\}_{n=0}^{\infty }$, $\left\{ {b_{n}}\right\}_{n=0}^{\infty }$ are given,
then we define  the canonical  polynomial system  by the recurrence relation
\begin{equation}
x {\psi_{n}}(x)= {b_{n}{\psi_{n+1}}(x)}+{a_{n}{\psi_{n}}(x)}+{b_{n-1}{\psi_{n-1}}(x)},\qquad
n\geq 0,\qquad
b_{-1}=0,
\label{trivs4}
\end{equation}

\begin{equation}
{\psi_{0}(x)}=1.
\label{trivs5}
\end{equation}
As before, the remark {\ref{j}} is true.
\begin{thm}
\label{gh}
Let  $\left\{ {\psi_{n}(x)}\right\}_{n=0}^{\infty }$  be a real
polynomials system defined by
(\ref{trivs4}), (\ref{trivs5}), and let $\mu$  be a probability
measure on $[a,b]\subset R$.
 A system of polynomials  $\left\{ {\psi_{n}(x)}\right\}_{n=0}^{\infty }$
is orthonormal
with respect to the measure  $\mu$ on $[a,b]\subset R$ if and only if
the coefficients
$\left\{ {a_{n}}\right\}_{n=0}^{\infty }$,
$\left\{ {b_{n}}\right\}_{n=0}^{\infty }$ involved
in the recurrence relations  (\ref{trivs4}) are the solution of the system
(\ref{trivs2}), (\ref{trivs3}),
where $\mu_{k}$ are defined by  (\ref{trivs1}).
\end{thm}
\begin{rem}
Just as in subsection {2.2} we introduce the Hilbert space ${\tt G}_{1}$ and the set
$\left\{ {\phi_{n}(x)}\right\}_{n=0}^{\infty }$. The corollary {\ref{ded}} still stands for this system.
\end{rem}

\subsection{}

We determine (just as in subsection {2.3}) the Poisson kernel in the Hilbert space
${\tt F}_{1}\otimes {\tt F}_{2}$ and the operators $U_x$  and  $U_y$. The lemmas {\ref{c}} and {\ref{u}}
still stand for these operators. As before, one can
define  the momentum operator $P_{{\tt F}_{1}}$,
which is conjugate to the position operator $X_{{\tt F}_{1}}$ with respect to the basis
$\left\{ {\varphi_{n}(x)}\right\}_{n=0}^{\infty }$ in  ${\tt F}_{1}$, and the symmetric hamiltonian ${H_{{\tt F}_{1}}}(t)$,
 which does not have to be a selfadjoint operator. Moreover, the set  $\left\{ {\varphi_{n}(x)}\right\}_{n=0}^{\infty }$
does not have to be a set of eigenfunctions of the operator  ${H_{{\tt F}_{1}}}(t)$ at any value  $t$. However one can
 remedy the situation by using  new position and momentum operators.

As in  the lemma {\ref{du}} we have
\begin{align}
\label{trivs6}
X_{{\tt F}_{1}}{\varphi_{0}(x)}&={b_{0}}{\varphi_{1}(x)}+{a_{0}}{\varphi_{0}(x)},\\
\label{trivs7}
P_{{\tt F}_{1}}{\varphi_{0}(x)}&=-(\imath){b_{0}}{\varphi_{1}(x)}+{a_{0}}{\varphi_{0}(x)},\\
\label{trivs8}
X_{{\tt F}_{1}}{\varphi_{n}(x)}&={b_{n-1}}{\varphi_{n-1}(x)}+{a_{n}}{\varphi_{n}(x)}+{b_{n}}{\varphi_{n+1}(x)},\quad
n\geq 1,\\
\label{trivs9}
P_{{\tt F}_{1}}{\varphi_{n}(x)}&={\imath}({b_{n-1}\varphi_{n-1}(x)}-{{b_{n}\varphi_{n+1}(x)}})+{a_{n}}{\varphi_{n}(x)},
n\geq 1.
\end{align}
From (\ref{trivs6}) -  (\ref{trivs9})  it follows that:
\begin{align}
\label{trivs10}
(X_{{\tt F}_{1}}-P_{{\tt H}_{1}}){\psi_{0}(x)}&={b_{0}}{\psi_{1}(x)},\\
\label{trivs11}
(X_{{\tt F}_{1}}-P_{{\tt H}_{1}}){\psi_{n}(x)}&=
(1-\imath){b_{n-1}}{\psi_{n-1}(x)}+(1+\imath){b_{n}}{\psi_{n+1}(x)},\quad
n\geq 1.
\end{align}
When (\ref{trivs10}), (\ref{trivs11}) is compared with (\ref{triv25}),(\ref{triv26}),(\ref{triv28}),(\ref{triv29}),
it is apparent that one can  introduce the new position $\widetilde{X}_{{\tt F}_{1}}$
and momentum $\widetilde{P}_{{\tt F}_{1}}$ operators as follows:
\begin{align}
\label{trivs12}
\widetilde{X}_{{\tt F}_{1}}&=Re(X_{{\tt F}_{1}}-P_{{\tt F}_{1}}),\\
\label{trivs13}
\widetilde{P}_{{\tt F}_{1}}&=(-\imath)Im(X_{{\tt F}_{1}}-P_{{\tt F}_{1}}).
\end{align}
If we replace $X_{{\tt F}_{1}}\longmapsto \widetilde{X_{{\tt F}_{1}}}$
and  ${P_{{\tt F}_{1}}\longmapsto {\widetilde{P}_{{\tt F}_{1}}}}$ ,
then the formulas (\ref{triv25}),(\ref{triv26}),(\ref{triv28}),
(\ref{triv29}) are valid for the operators  $\widetilde{X}_{{\tt F}_{1}}$ and  $\widetilde{P}_{{\tt F}_{1}}$.
\begin{lem}
Let the operators $\widetilde{X}_{{\tt H}_{1}}$ and $\widetilde{P}_{{\tt H}_{1}}$ be defined by
(\ref{trivs12}), (\ref{trivs13}). Then we have the formula  (\ref{triv20}) (with $t=\imath$)
\begin{equation}
\widetilde{P}_{{\tt F}_{1}}=K_{\tt F}^{*}{\widetilde{Y}_{{\tt F}_{2}}}{K_{\tt F}}.
\label{trivs14}
\end{equation}
\end{lem}
Now we define the energy operator:
\begin{equation}
\widetilde{H}_{{\tt F}_{1}}={{\widetilde{X}_{{\tt F}_{1}}}^{2}}+{{\widetilde{P}_{{\tt F}_{1}}}^{2}}.
\label{trivs15}
\end{equation}
The following theorem is similar to  theorem {\ref{b}}.
\begin{thm}
The  operator $\widetilde{H}_{{\tt F}_{1}}$ defined by (\ref{trivs15}) is a selfadjoint operator in the
space ${\tt F}_{1}$ with a orthonormal basis  $\left\{ {\varphi_{n}(x)}\right\}_{n=0}^{\infty }$. Moreover
the set $\left\{ {\varphi_{n}(x)}\right\}_{n=0}^{\infty }$ is a eigenfunction system of the operator
$\widetilde{H}_{{\tt F}_{1}}$ and the eigenvalues of this operator are equal to:
\begin{equation}
\lambda_{0}=2{b_{0}^{2}},\qquad
\lambda_{n}=2({b_{n-1}^{2}}+{b_{n}^{2}}),\qquad
n\geq 1.
\label{trivs16}
\end{equation}
\end{thm}
We define the ladder operators:
\begin{equation}
\widetilde{a^{+}_{{\tt H}_{1}}}={\frac{1}{\sqrt{2}}}
\left (
{\widetilde{X}_{{\tt H}_{1}}}+\imath{\widetilde{P}_{{\tt H}_{1}}}
\right )
,\qquad
\widetilde{a^{-}_{{\tt H}_{1}}}={\frac{1}{\sqrt{2}}}
\left (
{\widetilde{X}_{{\tt H}_{1}}}-\imath{\widetilde{P}_{{\tt H}_{1}}}
\right ),
\label{trivs17}
\end{equation}
If we replace  $a^{\pm}_{{\tt H}_{1}}\longmapsto \widetilde{a^{\pm}_{{\tt H}_{1}}}$),  then the formulas
(\ref{triv36}) are valid. Moreover  the theorem {\ref{g}} is also true.
\begin{rem}
It should be stressed that in this case, too, we succeeded in constructing  some oscillator system. However, now the
position operator does not have to be an operator of the multiplication on an independent variable.
\end{rem}

\subsection{}

Now we consider a nonsymmetric Jacobi matrix of a position operator in Fock representation. Let sequences
$\left\{ {a_{n}}\right\}_{n=0}^{\infty }$, $\left\{ {b_{n}}\right\}_{n=0}^{\infty }$,$\left\{ {c_{n}}\right\}_{n=0}^{\infty }$
and a sequence $\left\{ {A_{k,n}}\right\}_{k,n=0}^{\infty }$ be a solution to the following equation system:
\begin{equation}
A_{k,n}= {b_{n}{A_{k-1,n+1}}}+{a_{n}{A_{k-1,n}}}+{c_{n}{A_{k-1,n-1}}},\qquad
n\geq 0,\qquad
b_{-1}=0,
\label{trivs18}
\end{equation}
satisfying the initial conditions  (\ref{trivs3}) too.
Contrary to  (\ref{trivs2}) there is an infinite family of solution to the system (\ref{trivs18}) , (\ref{trivs3}).
We can find uniquely from (\ref{trivs18}) , (\ref{trivs3}) only:
\begin{equation}
d_{j}=\sqrt{{b _{j}}{c_{j+1}}},\qquad
j= 0,1,\dots .
\end{equation}
If the sequences $\left\{ {a_{n}}\right\}_{n=0}^{\infty }$, $\left\{ {b_{n}}\right\}_{n=0}^{\infty }$,$\left\{ {c_{n}}\right\}_{n=0}^{\infty }$
are given, then we define the  polynomial system
$\left\{ {\hat{\psi}_{n}}(x)\right\}_{n=0}^{\infty }$ by:
\begin{equation}
x {\hat{\psi}_{n}}(x)= {b_{n}{\hat{\psi}_{n+1}}(x)}+{a_{n}{\hat{\psi}_{n}}(x)}+{c_{n}{\hat{\psi}_{n-1}}(x)},\qquad
n\geq 0,\qquad
c_{0}=0,
\label{trivs19}
\end{equation}

\begin{equation}
{\hat{\psi}_{0}(x)}=1.
\label{trivs20}
\end{equation}
If the sequences $\left\{ {a_{n}}\right\}_{n=0}^{\infty }$, $\left\{ {d_{n}}\right\}_{n=0}^{\infty }$,
are given, then  the canonical polynomial system $\left\{ {\psi_{n}}(x)\right\}_{n=0}^{\infty }$ is
defined  by the recurrence relations (\ref{trivs4}),(\ref{trivs5})  with  $d_{n}$    instead of  $b_{n}$.
It follows from the theorem {\ref{gh}} that the set $\left\{ {\psi_{n}}(x)\right\}_{n=0}^{\infty }$ is
an orthonormal polynomials system in the space ${\tt H}_{1}$ .
It can easily be checked that the renormalizaton:
\begin{equation}
{\hat{\psi}_{n}}(x)={\gamma _{n}}{\psi_{n}}(x),\qquad
n\geq 0,
\label{trivs21}
\end{equation}
where
\begin{equation}
{\gamma _{0}}=1,\quad
{\gamma _{n}}=\sqrt{\frac{c_{1}\cdot c_{2}\cdots c_{n} }{b_{0}\cdot b_{1}\cdots b_{n-1}}},\quad
n\geq 1,
\label{trivs22}
\end{equation}
reduce  (\ref{trivs19}), (\ref{trivs20}) to the symmetric relations (\ref{trivs4}),(\ref{trivs5}).
From the orthonormal conditions for the system $\left\{ {\psi_{n}}(x)\right\}_{n=0}^{\infty }$
we obtain the following orthogonal relations:
\begin{equation}
\int_{-\infty }^{\infty }{\hat{\psi}_{n}(x){\hat{\psi}_{m}(x)}{\mu (dx)}}={\gamma _{n}}{\gamma _{m}}
\int_{-\infty }^{\infty }{{\psi_{n}(x)}{{\psi_{m}(x)}}{\mu (dx)}}=
{\gamma _{n}}^{2}{\delta _{nm}}.
\label{trivs23}
\end{equation}
Note that the remark {\ref{gg}} is true in this case  too.
A main example of the nonsymmetric scheme for the classical orthogonal polynomials
is the Laguerre polynomials.

\section{Laguerre polynomials}

Denote by  ${{\tt G}_{1}}={L^{2}}({R^{1}}_{+}),$\quad  ${{\tt H}_{1}}={L^{2}}({R^{1}}_{+};{x^{\alpha }}{\exp(-x)}{dx}$
and
\begin{equation}
\psi (x)= {x ^{\frac{\alpha }{2}}}\exp(-\frac{x}{2}).
\label{trivs24}
\end{equation}
We determine the Laguerre polynomials ${{L^{\alpha }}_{n}}(x)$  (\cite{sze},\cite{koe}):
\begin{equation}
{L^{\alpha }_{n}}(x)={\frac{(\alpha +1)_{n}}{n!}}{}_{1}F_{1}
\left (
-n;\alpha +1;x
\right ).
\label{trivs25}
\end{equation}
Let
\begin{equation}
{d_{n}}^{2}=\frac{\Gamma (n+\alpha +1)}{n!}.
\label{trives26}
\end{equation}
We define also the orthonormal systems $\left\{ {\psi_{n}(x)}\right\}_{n=0}^{\infty }$
and  $\left\{ {\phi_{n}(x)}\right\}_{n=0}^{\infty }$ by the following  formulas:
\begin{equation}
{\psi_{n}}(x)= {{d_{n}}^{-1}}{L^{\alpha }_{n}}(x),\qquad
{\phi_{n}}(x)=\psi(x){\psi_{n}}(x),\qquad
n\geq 0.
\label{trivs26}
\end{equation}
Using the recurrence relations for the Laguerre polynomials
\begin{align}
(n+1){L^{\alpha }_{n+1}}(x)&=
(2n+\alpha  +1-x){{L^{\alpha }}_{n}}(x)+
(n+\alpha ){L^{\alpha }_{n-1}}(x),\quad
n\geq 1,\\
{L^{\alpha }_{0}}(x)&=1,\qquad
{L^{\alpha }_{-1}}(x)=0,
\end{align}
we have  (\ref{trivs19}), where
\begin{align}
b_{n}&=-(n+1)\frac{d_{n+1}}{d_{n}},\notag\\
a_{n}&=2n+\alpha +1,\notag\\
c_{n}&=-(n+\alpha )\frac{d_{n-1}}{d_{n}}.
\label{trivs27}
\end{align}
Finally we obtain from  (\ref{trivs27}) and  (\ref{trivs26}):
\begin{align}
b_{n}&=-\sqrt{(n+1)(n+\alpha +1)} ,\notag\\
a_{n}&=2n+\alpha +1,\notag\\
c_{n}&=b_{n-1}.
\label{trivs28}
\end{align}
We consider a differential operator $K$, which shall play a large role below. The operator $K$
acts on basis vectors by the following formulas (\cite{sze}):
 \begin{equation}
K={x}\frac{d}{dx}= {b_{n-1}}{\psi_{n-1}}(x)+n{\psi_{n}}(x),\qquad
n\geq 0,\quad
b_{-1}=0.
\label{trivs29}
\end{equation}
From  (\ref{trivs29}) and  (\ref{trivs10}) , (\ref{trivs11})
we get the formula for the operator  $P_{{\tt H}_{1}}$:
\begin{equation}
P_{{\tt H}_{1}}=2N_{{\tt H}_{1}}+(\alpha +1)I_{{\tt H}_{1}}-
\imath(2K-X_{{\tt H}_{1}}+(\alpha +1)I_{{\tt H}_{1}})
\label{trivs30}
\end{equation}
Combining  (\ref{trivs30}) , (\ref{trivs12}) and (\ref{trivs13}) we have:
\begin{align}
\label{trivs31}
\widetilde{X}_{{\tt H}_{1}}&=X_{{\tt H}_{1}}-2N_{{\tt H}_{1}}-(\alpha +1)I_{{\tt H}_{1}},\\
\label{trivs32}
\widetilde{P}_{{\tt H}_{1}}&=\imath(2K-X_{{\tt H}_{1}}+(\alpha +1)I_{{\tt H}_{1}}).
\end{align}
Further, from (\ref{trivs29}) and the definition (\ref{triv35}) and (\ref{triv38}) of the operators
 $a^{\pm}_{{\tt H}_{1}}$, $N_{{\tt H}_{1}}$ it follows that:
\begin{align}
\label{trivs33}
\widetilde{a^{-}_{{\tt H}_{1}}}&=\sqrt{2} (K-N_{{\tt H}_{1}}),\\
\label{trivs34}
\widetilde{a^{+}_{{\tt H}_{1}}}&=\sqrt{2} ([N_{{\tt H}_{1}},X_{{\tt H}_{1}}]+(K-N_{{\tt H}_{1}})).
\end{align}
In view of
\begin{align}
\label{trivs35}
[K,X_{{\tt H}_{1}}]&=X_{{\tt H}_{1}},\\
\label{trivs36}
[K-N_{{\tt H}_{1}},X_{{\tt H}_{1}}]&=2K+(\alpha +1)I_{{\tt H}_{1}}.
\end{align}
it is not hard to prove that (\ref{trivs17}) is true for the operators (\ref{trivs31})-(\ref{trivs34}).
Then one can rewrite (\ref{trivs31}),(\ref{trivs32}) in the form:
\begin{align}
\label{trivs37}
\widetilde{X}_{{\tt H}_{1}}&=2(K-N_{{\tt H}_{1}})+[N_{{\tt H}_{1}},X_{{\tt H}_{1}}],\\
\label{trivs38}
\widetilde{P}_{{\tt H}_{1}}&=(-\imath)[N_{{\tt H}_{1}},X_{{\tt H}_{1}}].
\end{align}
Taking into account (\ref{trivs37}),(\ref{trivs38}) and (\ref{trivs33}),(\ref{trivs34})
one can write the hamiltonian $\widetilde{H}_{{\tt H}_{1}}$ defined by  (\ref{trivs15}) in the following form:
\begin{align}
\widetilde{H}_{{\tt H}_{1}}&=4(K-N_{{\tt H}_{1}})^{2}+
2([N_{{\tt H}_{1}},X_{{\tt H}_{1}}](K-N_{{\tt H}_{1}})+\notag\\
+(K-N_{{\tt H}_{1}})[N_{{\tt H}_{1}},X_{{\tt H}_{1}}])&=
\widetilde{a^{+}_{{\tt H}_{1}}}\widetilde{a^{-}_{{\tt H}_{1}}}+
\widetilde{a^{-}_{{\tt H}_{1}}}\widetilde{a^{+}_{{\tt H}_{1}}}.
\label{trivs39}
\end{align}
Using  (\ref{trivs30}),(\ref{trivs31}) we also have:
\begin{align}
\widetilde{H}_{{\tt H}_{1}}=(X_{{\tt H}_{1}}&-2N_{{\tt H}_{1}})^{2}-
2(\alpha +1)(X_{{\tt H}_{1}}-2N_{{\tt H}_{1}})\notag\\
-(2K-X_{{\tt H}_{1}})^{2}&-2(\alpha +1)(2K-X_{{\tt H}_{1}})=4({N_{{\tt H}_{1}}}^{2}-K^{2})\notag\\
+2((K-N_{{\tt H}_{1}})X_{{\tt H}_{1}}&+X_{{\tt H}_{1}}(K-N_{{\tt H}_{1}}))+
4(\alpha +1)(N_{{\tt H}_{1}}-K).
\label{trivs40}
\end{align}
Moreover, the energy levels are
\begin{equation}
\lambda _{n}=4(n^{2}+(\alpha +1)n+\frac{\alpha +1}{2}),\qquad
n\geq 0.
\label{trivs41}
\end{equation}
The following theorem is valid.
\begin{thm}
\label{dor}
The equation
\begin{equation}
{\widetilde{H}_{{\tt H}_{1}}}{\psi_{n}}(x)={\lambda _{n}}{\psi_{n}}(x),\qquad
n\geq 0,
\label{trivs42}
\end{equation}
where  $\lambda _{n}$ is defined by  (\ref{trivs41}), is equivalent to the differential equation
for the Laguerre polynomials:
\begin{equation}
x({L_{n}^{(\alpha )}}(x))^{\prime \prime }+
(\alpha +1-x)({L_{n}^{(\alpha )}}(x))^{\prime }+n{{L{n}^{(\alpha )}}(x)}=0,\quad
n\geq 0.
\label{trivs43}
\end{equation}
\end{thm}
\begin{proof}
Using (\ref{trivs29}) and (\ref{trivs26}) we rewrite the differential equation (\ref{trivs43})
in the form of the operator equality in the space ${\tt H}_{1}$:
\begin{equation}
K^{2}+{\alpha }K-{X_{{\tt H}_{1}}}(K-{N_{{\tt H}_{1}}})=0.
\label{trivs44}
\end{equation}
In view of (\ref{trivs36}) and  (\ref{trivs37}) the equation (\ref{trivs42})  is equivalent to the following
operator equality in the space ${\tt H}_{1}$:
\begin{equation}
K^{2}+{\alpha }K-{X_{{\tt H}_{1}}}(K-{N_{{\tt H}_{1}}})=\frac{1}{2}[K-N_{{\tt H}_{1}},X_{{\tt H}_{1}}]
-K-\frac{\alpha +1}{2}{I_{{\tt H}_{1}}}.
\label{trivs45}
\end{equation}
It is obvious from (\ref{trivs44}) and (\ref{trivs45}) that it is sufficient to prove that  the right-hand
side of  (\ref{trivs45}) vanishes.From (\ref{trivs36}) it follows that the latter is true.
\end{proof}
\begin{rem}
The theorem {\ref{g}} is true in our case. Then we have
\begin{equation}
\label{trivs46}
[\widetilde{a^{-}_{{\tt H}_{1}}},\widetilde{a^{+}_{{\tt H}_{1}}}]{\psi_{n}}(x)=
2({b_{n}}^{2}-{b_{n-1}}^{2}){\psi_{n}}(x),\quad
n\geq 0,\quad
b_{-1}=0.
\end{equation}
Taking into account (\ref{trivs28}) we calculate
\begin{equation}
\label{trivs47}
{b_{n}}^{2}-{b_{n-1}}^{2}=2n+\alpha +1.
\end{equation}
From (\ref{trivs47}) and  (\ref{trivs46}) we get  the following commutation relation:
\begin{equation}
\label{trivs48}
[\widetilde{a^{-}_{{\tt H}_{1}}},\widetilde{a^{+}_{{\tt H}_{1}}}]=
2{N_{{\tt H}_{1}}}+(\alpha +1){I_{{\tt H}_{1}}}.
\end{equation}
\end{rem}
As another instance of the nonsymmetric scheme we consider the Jacobi polynomials
${P_{n}^{(\alpha ,\beta )}}(x)$   under the condition
$\alpha\not=  \beta $.

\section{The Jacobi polynomials}

The Jacobi polynomials   (\cite{sze}) one can be determined by
\begin{equation}
{P_{n}^{(\alpha ,\beta )}}(x)={\frac{(\alpha +1)_{n}}{n!}}
{{}_{2}F_{1}
\left (
-n,n+\alpha +\beta +1;\alpha +1;\frac{1-x}{2}
\right )}.
\label{trivs49}
\end{equation}
Let
$$
{{\tt G}_{1}}={L^{2}}([-1,1]),
$$
$$
{{\tt H}_{1}}={L^{2}}([-1,1];{({d_{0}}(\alpha ,\beta ))^{-2}}{(1-x)^{\alpha }}{(1+x)^{\beta }}dx),
$$
where
$$
{d_{0}^{2}}(\alpha ,\beta )=2^{{\alpha +\beta +1}}\frac{{\Gamma (\alpha +1)}{\Gamma (\beta  +1)}}{\Gamma (\alpha +\beta +2)}.
$$
 Let
\begin{equation}
\psi (x)={(1-x)^{\frac{\alpha }{2}}}{(1+x)^{\frac{\beta }{2}}}.
\end{equation}
We define the orthonormal systems  $\left\{ {\psi_{n}(x)}\right\}_{n=0}^{\infty }$ and $\left\{ {\phi_{n}(x)}\right\}_{n=0}^{\infty }$
by the following formulas:
\begin{equation}
{\psi_{n}}(x)={d_{0}}{d_{n}^{-1}}{P_{n}^{(\alpha ,\beta )}}(x),\qquad
{\phi_{n}}(x)={\psi(x)}{\psi_{n}}(x),\qquad
n\geq 0,
\label{trivs50}
\end{equation}
where the constants  $d_{n}$ are given by
\begin{equation}
d_{n}^{2}=\frac{{2^{\alpha +\beta +1}}{(\Gamma (n+\alpha +1))}{(\Gamma (n+\beta +1))}}
{(2n+\alpha \beta +1){n!}\Gamma (n+\alpha +\beta +1)},\quad
n\geq 0.
\label{trivs51}
\end{equation}
Using (\ref{trivs51}) and the recurrence relations for the Jacobi polynomials  (see \cite{sze})
we get  (\ref{trivs4}), (\ref{trivs5}), where
\begin{equation}
a_{n}={\frac{{\beta }^{2}-{\alpha }^{2}}{(2n+\alpha +\beta )(2n+2+\alpha +\beta )}},\quad
n\geq 0,
\label{trivs52}
\end{equation}
\begin{equation}
b_{n}=\sqrt{\frac{(n+1)(n+\alpha +1)(n+\beta +1)(n+\alpha +\beta +1)}
{(2n+\alpha +\beta +1)(2n+\alpha +\beta +2)^{2}(2n+\alpha +\beta +3)}},\quad
n\geq 0.
\label{trivs53}
\end{equation}
It is known how the operator $A$ defined by (\ref{triv68}) acts on the Jacobi polynomials (\cite{sze}).
Then it is not hard to get from (\ref{trivs50}),(\ref{trivs51}) the following equalities:
\begin{align}
A{\psi_{n}}(x)&=(n+\alpha +\beta +1){b_{n-1}}{\psi_{n-1}}(x)-\notag\\
\frac{2n(n+\alpha +\beta +1)}{\alpha +\beta }a_{n}{\psi_{n-1}}(x)&-
{n{b_{n}}{\psi_{n+1}}(x)},\qquad
n\geq 0,\quad
b_{-1}=0.
\label{trivs54}
\end{align}
Multiplying both sides of (\ref{trivs4}) by $n+\alpha +\beta +1$ and subtracting (\ref{trivs54})
from the obtained result we eliminate $\psi_{n-1}$ from (\ref{trivs54}). Then we obtain
\begin{align}
\frac{\sqrt{2} }{2n+\alpha +\beta +1}
(x(n+\alpha +\beta +1)&-A-\frac{a_{n}(n+\alpha +\beta +1)(2n+\alpha +\beta )}{\alpha +\beta })\notag\\
{\psi_{n}}(x)&={\sqrt{2} }b_{n}{\psi_{n+1}}(x),\qquad
n\geq 0.
\label{trivs55}
\end{align}
From  (\ref{trivs55}),(\ref{trivs17}) and (\ref{triv36}) it follows that
\begin{align}
\label{trivs56}
\widetilde{a^{+}_{{\tt H}_{1}}}&=\sqrt{2} ((X_{{\tt H}_{1}}(N_{{\tt H}_{1}}+(\alpha +\beta +1)I_{{\tt H}_{1}})-A)\notag\\
(2N_{{\tt H}_{1}}&+(\alpha +\beta +1)I_{{\tt H}_{1}})^{-1}+(\alpha -\beta )
(N_{{\tt H}_{1}}+(\alpha +\beta +1)I_{{\tt H}_{1}})\notag\\
(2N_{{\tt H}_{1}}&+(\alpha +\beta +1)I_{{\tt H}_{1}})^{-1} (2N_{{\tt H}_{1}}+(\alpha +\beta +2)I_{{\tt H}_{1}})^{-1}).
\end{align}
In order to eliminate  $\psi_{n+1}$ from (\ref{trivs54}) we  multiply both sides of (\ref{trivs4}) by $n$ and
add this to (\ref{trivs54}). Then we get
\begin{align}
\frac{\sqrt{2} }{2n+\alpha +\beta +1}
(xn&+A+\frac{a_{n}n(2n+\alpha +\beta +2)}{\alpha +\beta })\notag\\
{\psi_{n}}(x)&={\sqrt{2} }b_{n-1}{\psi_{n-1}}(x),\qquad
n\geq 0.
\label{trivs57}
\end{align}
Combining  (\ref{trivs57}),(\ref{trivs17}) and  (\ref{triv36}) we have
\begin{align}
\label{trivs58}
\widetilde{a^{-}_{{\tt H}_{1}}}&=\sqrt{2} (X_{{\tt H}_{1}}N_{{\tt H}_{1}}+A-(\alpha -\beta )N_{{\tt H}_{1}}\notag\\
(2N_{{\tt H}_{1}}&+(\alpha +\beta )I_{{\tt H}_{1}})^{-1})
(2N_{{\tt H}_{1}}+(\alpha +\beta +1)I_{{\tt H}_{1}})^{-1}.
\end{align}
Taking into account (\ref{trivs17})  and  (\ref{trivs56}),(\ref{trivs58}) one can write
\begin{align}
\label{trivs59}
\widetilde{X}_{{\tt H}_{1}}&=\frac{1}{\sqrt{2} }(\widetilde{a^{+}_{{\tt H}_{1}}}+
\widetilde{a^{-}_{{\tt H}_{1}}})=X_{{\tt H}_{1}}N_{{\tt H}_{1}}+\notag\\
(\alpha -\beta )(2N_{{\tt H}_{1}}&+(\alpha +\beta +1)I_{{\tt H}_{1}})^{-1}
((2\alpha +2\beta +1)N_{{\tt H}_{1}}+(\alpha +\beta +1)^{2})\notag\\
(2N_{{\tt H}_{1}}&+(\alpha +\beta )I_{{\tt H}_{1}})^{-1})
(2N_{{\tt H}_{1}}+(\alpha +\beta +2)I_{{\tt H}_{1}})^{-1}),\\
\label{trivs60}
\widetilde{P}_{{\tt H}_{1}}&=\frac{-\imath}{\sqrt{2} }(\widetilde{a^{+}_{{\tt H}_{1}}}-
\widetilde{a^{-}_{{\tt H}_{1}}})={-\imath}(-2A+X_{{\tt H}_{1}}(\alpha +\beta +1))+\notag\\
(\alpha -\beta )((2N_{{\tt H}_{1}}&+(\alpha +\beta +1)I_{{\tt H}_{1}})^{2}+N_{{\tt H}_{1}})
(2N_{{\tt H}_{1}}+(\alpha +\beta )I_{{\tt H}_{1}})^{-1}\notag\\
(2N_{{\tt H}_{1}}&+(\alpha +\beta +2)I_{{\tt H}_{1}})^{-1})
(2N_{{\tt H}_{1}}+(\alpha +\beta +1)I_{{\tt H}_{1}})^{-1}.
\end{align}
It is not hard to find from (\ref{trivs59}),(\ref{trivs60}) and the definition  (\ref{trivs15})
the explicit form of the hamiltonian  $\widetilde{H}_{{\tt H}_{1}}$. In view of
(\ref{trivs16}),(\ref{trivs53}) we get the energy levels by
\begin{align}
\label{trivs61}
\lambda _{n}&=2(b_{n-1}^{2}+b_{n}^{2})=\notag\\
{}&=\frac{(2n+\alpha +\beta +1)^{2}(s_{n}-4w_{n})+5s_{n}-2w_{n}}
{(2n+\alpha +\beta )^{2}(2n+\alpha +\beta +2)^{2}(2n+\alpha +\beta -1)(2n+\alpha +\beta +3)},
\end{align}
where
\begin{align}
\label{trivs62}
s_{n}&=t_{n}+t_{n+1},\qquad
t_{n}=n(n+\alpha )(n+\beta )(n+\alpha +\beta ),\notag\\
w_{n}&=2n^{2}+2n(\alpha +\beta +1)+(\alpha +1)(\beta +1).
\end{align}
The next theorem is an extension of the theorem {\ref{dof}}.
\begin{thm}
\label{dob}
The equation $\widetilde{H}_{{\tt H}_{1}}{\psi_{n}}(x)={\lambda _{n}}{\psi_{n}}(x),\quad n\geq 0$,
where $\lambda _{n}$ defined by  (\ref{trivs61}) is equivalent to the differential equation
for the Jacobi polynomials:
\begin{align}
\label{trivs63}
\frac{d}{dx}(({1}&-x)^{\alpha +1}(1+x)^{\beta +1}\frac{d}{dx}){P_{n}^{(\alpha ,\beta  )}}(x)+\notag\\
+n(n+\alpha &+\beta +1)(1-x)^{\alpha }(1+x)^{\beta }{P_{n}^{(\alpha ,\beta )}}(x)=0,
\end{align}
$(n\geq 0.)$
\end{thm}
\begin{rem}
The result may be proved in much the same way as the theorem {\ref{dof}}. Here we omit this
proof as well as the explicit forms for the number operator $\widetilde{N}_{{\tt H}_{1}}$, the momentum
operator $\widetilde{P}_{{\tt H}_{1}}$ and the hamiltonian $\widetilde{H}_{{\tt H}_{1}}$.
\end{rem}
In conclusion we will point out some associations  between the canonical systems in the symmetric
and nonsymmetric schemes.

\section{Connection of symmetric with nonsymmetric schemes}

Denote by ${\mu }^{s}$ - a symmetric probability measure on $R^{1}$ and by
${\tt H}^{s}_{1}=L^{2}(R^{1};{\mu }^{s})$. Let  ${\mu }^{s}(0)=0$ and
\begin{equation}
\label{trivs64}
{\mu }^{s}={\mu }_{+}+{\mu }_{-}
\end{equation}
be a decomposition of the measure ${\mu }^{s}$ into the orthogonal sum of two (nonsymmetric) measures
${\mu }_{+}$ and ${\mu }_{-}$ defined by the equalities:
\begin{equation}
\label{trivs65}
{\mu }_{+}(B)={\mu }^{s}(R^{1}_{+}\cap B),\qquad
{\mu }_{-}(B)={\mu }^{s}(R^{1}_{-}\cap B),
\end{equation}
for any Borel set  $B\subset R^{1}$.
Let $\left\{ {\psi_{n}(x)}\right\}_{n=0}^{\infty }$ be a canonical (complete) orthonormal polynomial
system in the space ${\tt H}^{s}_{1}$. Suppose that  this system is constructed by the sequence
$\left\{ b^{s}_{n}\right\}_{n=0}^{\infty}$ via the recurrence relations (\ref{triv5}) ,(\ref{triv6}).
Further, we denote by
\begin{equation}
\label{trivs66}
{\tt H}^{+}_{1}=L^{2}(R^{1}_{+};2{\mu }_{+}(dx)),\quad
{\tt H}^{-}_{1}=L^{2}(R^{1}_{-};2{\mu }_{-}(dx)).
\end{equation}
Let $\left\{ {{\psi }^{+}_{n}(x)}\right\}_{n=0}^{\infty }$, $\left\{ {{\psi }^{-}_{n}(x)}\right\}_{n=0}^{\infty }$
be canonical (complete) orthonormal polynomial systems  in ${\tt H}^{+}_{1}$ and ${\tt H}^{-}_{1}$ respectively.
Let  these systems be constructed by  real sequences $\left\{ a^{a}_{n}\right\}_{n=0}^{\infty}$ ,
 $\left\{ b^{a}_{n}\right\}_{n=0}^{\infty}$ via the recurrence relations:
\begin{equation}
x{\psi }^{+}_{n}(x)=b^{a}_{n}{\psi }^{+}_{n+1}(x)+a^{a}_{n}{\psi }^{+}_{n}(x)+b^{a}_{n-1}{\psi }^{+}_{n-1}(x),\quad
n\geq 0,\quad
b^{a}_{-1}=0,
\label{trivs67}
\end{equation}
\begin{equation}
{\psi }^{+}_{0}(x)=1,
\label{trivs68}
\end{equation}
and
\begin{equation}
x{\psi }^{-}_{n}(x)=b^{a}_{n}{\psi }^{-}_{n+1}(x)-a^{a}_{n}{\psi }^{-}_{n}(x)+b^{a}_{n-1}{\psi }^{-}_{n-1}(x),\quad
n\geq 0,\quad
b^{a}_{-1}=0,
\label{trivs69}
\end{equation}
\begin{equation}
{\psi }^{-}_{0}(x)=1.
\label{trivs70}
\end{equation}
The following lemmas are valid. The proof of these is left to the reader.
\begin{lem}
\label{dtf}
Let $2{\mu }_{+}$ and $2{\mu }_{-}$  be probability measures on $R^{1}_{+}$
and $R^{1}_{-}$ respectively such that the measure ${\mu }^{s}={\mu }_{+}+{\mu }_{-}$
is a symmetric probability measure on $R^{1}$. Let   $\left\{ {{\psi }^{+}_{n}(x)}\right\}_{n=0}^{\infty }$ and
$\left\{ {{\psi }^{-}_{n}(x)}\right\}_{n=0}^{\infty }$ be the orthonormal polynomial systems constructed by
$2{\mu }_{+}$ and $2{\mu }_{-}$ respectively via the recurrence relations (\ref{trivs67}) -(\ref{trivs70}).
Then the polynomial system $\left\{ {\psi_{n}(x)}\right\}_{n=0}^{\infty }$
defined by:
\begin{equation}
\label{trivs71}
\psi_{n}(x)=\frac{1}{2}({\psi }^{+}_{n}(x)+{\psi }^{-}_{n}(x)),\quad
n\geq 0,
\end{equation}
is an orthonormal system in the space
\begin{equation}
\label{trivs72}
{\tt H}^{s}_{1}={\tt H}^{+}_{1}\oplus {\tt H}^{-}_{1},
\end {equation}
where the measure ${\mu }^{s}$ is determined by (\ref{trivs64}) and satisfies the recurrence relations
 (\ref{triv5}) ,(\ref{triv6}) with the coefficients
$\left\{ b^{a}_{n}\right\}_{n=0}^{\infty}$.
\end{lem}
\begin{lem}
Let $\left\{ {\psi_{n}(x)}\right\}_{n=0}^{\infty }$ be a canonical orthonormal polynomial system in ${\tt H}^{s}_{1}$,
constructed by a symmetric probability measure ${\mu }^{s}$ (for more details we refer the reader to
the section {2}) via the recurrence relations (\ref{triv5}) ,(\ref{triv6}).
Then the polynomial systems  $\left\{ {{\psi }_{2l}(x)}\right\}_{l=0}^{\infty }$
and  $\left\{ {{\psi }_{2l+1}(x)}\right\}_{l=0}^{\infty }$
are the orthonormal  (it is understood that they are incomplete) systems in the spaces ${\tt H}^{+}_{1}$
and ${\tt H}^{-}_{1}$ respectively.
\end{lem}
\begin{rem}
1.One can obtain the corresponding complete orthonormal systems
 $\left\{ {{\psi }^{+}_{n}(x)}\right\}_{n=0}^{\infty }$ in  ${\tt H}^{+}_{1}$ and
$\left\{ {{\psi }^{-}_{n}(x)}\right\}_{n=0}^{\infty }$    in  ${\tt H}^{-}_{1}$ from
the set  $\left\{ {\psi_{n}(x)}\right\}_{n=0}^{\infty }$ by the Schmidt orthogonalization
in ${\tt H}^{+}_{1}$ and ${\tt H}^{-}_{1}$ respectively.

2. There is a simple relation between the moments  ${\mu }^{s}_{k}$
and ${\mu }^{+}_{k}$, ${\mu }^{-}_{k}$;
however connections between corresponding recurrence relations as well as between corresponding
oscillators are rather complicated.
\end{rem}

\section*{Acknowledgements}

This work has been supported by RFFI grant No  98-01-00310 .
The author is grateful to V.M. Babich, E.V. Damaskinsky, S.V. Kerov, I.V. Komarov,  and  P.P. Kulish
for discussions.

\bibliographystyle{amsplain}
\renewcommand{\refname}{References}

\bibliographystyle{amsplain}
\end{document}